\documentclass[12pt,reqno]{amsart}

\usepackage{amssymb}
\usepackage{bbm}
\usepackage{mathrsfs}
\usepackage[all]{xy}

\newtheorem{theo}{Theorem}[section]
\newtheorem{subl}[theo]{Sublemma}
\newtheorem{lem}[theo]{Lemma}
\newtheorem{prop}[theo]{Proposition}
\newtheorem{facs}[theo]{Facts}
\newtheorem{coro}[theo]{Corollary}

\theoremstyle{definition}
\newtheorem{defi}[theo]{Definition}
\newtheorem*{term}{Terminology}
\newtheorem*{mterm}{More terminology}
\newtheorem*{nota}{Notation}
\newtheorem*{conv}{Conventions and Notation}
\newtheorem*{tate}{The Tate conjecture for $K3$~surfaces}
\newtheorem*{mres}{The main result--A sufficient criterion for RM or CM}
\newtheorem*{link}{The link to the arithmetic}
\newtheorem*{fres}{Further results}
\newtheorem*{appl}{Applications}

\newtheorem*{proo1}{Proof of Theorem \ref{pc_yields_monodr}}
\newtheorem*{proo2}{Proof of Theorem \ref{Zarhin_rev}}
\newtheorem*{proo3}{Proof for Example \ref{Qw5_family}}

\theoremstyle{remark}
\newtheorem{ttt}[theo]{}
\newtheorem{rem}[theo]{Remark}
\newtheorem{rems}[theo]{Remarks}
\newtheorem{ex}[theo]{Example}
\newtheorem{exs}[theo]{Examples}

\newcommand{\br}{ }
\newcommand{\brr}{, }

\renewcommand{\atop}[2]{\genfrac{}{}{0pt}{}{#1}{#2}}

\newcommand{\pmodulo}[1]{\nobreak\mkern8mu
 (\textup{mod}\,\,#1)}

\newcommand{\notd}{\mathord{\nmid}}

\newcommand{\End}{\mathop{\text{\rm End}}\nolimits}
\newcommand{\Gal}{\mathop{\text{\rm Gal}}\nolimits}
\newcommand{\disc}{\mathop{\text{\rm disc}}\nolimits}

\newcommand{\Spec}{\mathop{\text{\rm Spec}}\nolimits}
\newcommand{\Frob}{\mathop{\text{\rm Frob}}\nolimits}
\newcommand{\Pic}{\mathop{\text{\rm Pic}}\nolimits}
\newcommand{\M}{\mathop{\text{\rm M}}\nolimits}
\newcommand{\GL}{\mathop{\text{\rm GL}}\nolimits}
\renewcommand{\O}{\mathop{\text{\rm O}}\nolimits}
\newcommand{\SO}{\mathop{\text{\rm SO}}\nolimits}
\newcommand{\GO}{\mathop{\text{\rm GO}}\nolimits}

\newcommand{\MG}{\mathop{\text{\rm MG}}\nolimits}
\newcommand{\Tr}{\mathop{\text{\rm Tr}}\nolimits}
\newcommand{\F}{\mathop{\text{\rm F}}\nolimits}
\newcommand{\rk}{\mathop{\text{\rm rk}}\nolimits}
\newcommand{\id}{\mathop{\text{\rm id}}\nolimits}
\newcommand{\im}{\mathop{\text{\rm im}}\nolimits}
\newcommand{\spann}{\mathop{\text{\rm span}}\nolimits}
\newcommand{\diag}{\mathop{\text{\rm diag}}\nolimits}

\newcommand{\alg}{\text{\rm alg}}
\newcommand{\tr}{\text{\rm tr}}
\newcommand{\pr}{\text{\rm pr}}
\newcommand{\et}{\text{\rm \'et}}
\newcommand{\Hg}{\text{\rm Hodge}}

\newcommand{\Pb}{{\text{\bf P}}}
\newcommand{\Ab}{{\text{\bf A}}}

\newcommand{\frakp}{{\mathfrak p}}

\newcommand{\frakA}{{\mathfrak A}}

\newcommand{\frakS}{{\mathfrak S}}
\newcommand{\frakX}{{\mathfrak X}}

\newcommand{\bbC}{{\mathbbm C}}
\newcommand{\bbF}{{\mathbbm F}}
\newcommand{\bbG}{{\mathbbm G}}
\newcommand{\bbN}{{\mathbbm N}}
\newcommand{\bbR}{{\mathbbm R}}
\newcommand{\bbQ}{{\mathbbm Q}}
\newcommand{\bbZ}{{\mathbbm Z}}

\newcommand{\calO}{{\mathscr{O}}}
\newcommand{\calP}{{\mathscr{P}}}
\newcommand{\calQ}{{\mathscr{Q}}}
\newcommand{\calR}{{\mathscr{R}}}
\newcommand{\calS}{{\mathscr{S}}}
\newcommand{\calT}{{\mathscr{T}}}
\newcommand{\calX}{{\mathscr{X}}}

\newcounter{abc}
\newenvironment{abc}{\begin{list}{\rm \alph{abc}) }%
{\usecounter{abc} \leftmargin=0.0pt \labelsep=0.0pt %
\listparindent=0.0pt \labelwidth=0.0pt \parsep=\smallskipamount%
 \itemsep=0.0pt \topsep=0.0pt \partopsep=\smallskipamount}}{\end{list}}

\newcounter{iii}
\newenvironment{iii}{\begin{list}{\rm \roman{iii}) }%
{\usecounter{iii} \leftmargin=0.0pt \labelsep=0.0pt %
\listparindent=0.0pt \labelwidth=0.0pt \parsep=\smallskipamount%
 \itemsep=0.0pt \topsep=0.0pt \partopsep=\smallskipamount}}{\end{list}}

\def\rightend#1#2{{%
 \leavevmode\nobreak\hskip .5em plus 1fil
 \penalty600 \hskip 0pt plus -1filll
 \vadjust{}\nobreak\hskip 0pt plus 1filll%
 #1\parfillskip=#2\relax \par}}

\def\eop{\ifmmode\rule[-22pt]{0pt}{1pt}\ifinner\tag*{$\square$}\else\eqno{\square}\fi\else\rightend{$\square$}{0pt}\fi}

\thanks{}

\title[Explicit families of
$K3$~surfaces
having real multiplication]{Explicit families of
{\boldmath$K3$}~surfaces
having \\ real multiplication}

\begin{document}

\author{Andreas-Stephan Elsenhans}

\address{Institut f\"ur Mathematik\\ Universit\"at W\"urzburg\\ Emil-Fischer-Stra\ss e 30\\ D-97074 W\"urzburg\\ Germany}
\email{stephan.elsenhans@mathematik.uni-wuerzburg.de}
\urladdr{https://www.mathematik.uni-wuerzburg.de/institut/personal/elsenhans.html}

\author[J\"org Jahnel]{J\"org Jahnel}

\address{\mbox{Department Mathematik\\ \!Univ.\ \!Siegen\\ \!Walter-Flex-Str.\ \!3\\ \!D-57068 \!Siegen\\ \!Germany}}
\email{jahnel@mathematik.uni-siegen.de}
\urladdr{http://www.uni-math.gwdg.de/jahnel}

\thanks{}

\date{January~23,~2020}

\keywords{$K3$~surface,
real multiplication, explicit examples}

\subjclass[2010]{14J28 primary; 14F20, 14F25, 14J10, 14J20, 11T06 secondary}

\begin{abstract}
For families of
$K3$
surfaces, we establish a sufficient criterion for real or complex multiplication. Our criterion is arithmetic in nature. It may show, at first, that the generic fibre of the family has a nontrivial endomorphism~field. Moreover, the endomorphism field does not shrink under specialisation. As an application, we present two explicit families of
$K3$
surfaces having real multiplication by
$\bbQ(\sqrt{2})$
and
$\bbQ(\sqrt{5})$,
respectively.
\end{abstract}

\maketitle
\thispagestyle{empty}

\section{Introduction}

Complex multiplication is a phenomenon that has been intensively studied, first and foremost for complex elliptic curves, cf.~\cite[Chapter~II]{Si} or~\cite[Chap\-ter~3]{Cox}. According~to its very definition, it is a purely geometric property. Nevertheless,~it has arithmetic consequences and much of the interest in complex multiplication stems from~these. The~generalisation to higher-dimensional abelian varieties is straightforward, except for the fact that, besides complex multiplication, the similar phenomena of real and quaternionic multiplication may occur. Abelian varieties are, however, not the~limit.

For~instance, let
$\frakX$
be a projective complex
$K3$~surface.
In~this situation, the occurrence of real and complex multiplication phenomena has been observed by Yu.\,G.\ Zarhin~\cite{Za}. They~are certainly deeper for
$K3$~surfaces
than for abelian varieties, because they do not concern the complex manifold directly, but merely its cohomology.

More~concretely, the cohomology
\mbox{$\bbQ$-vector}
space
$H^2(\frakX, \bbQ)$
is of
dimension~$22$.
On~the other hand, the rank of the Picard group
$\Pic \frakX$
may vary between
$1$
and~$20$.
Put~$P := c_1(\Pic \frakX) \!\otimes_\bbZ\! \bbQ \subset H^2(\frakX, \bbQ)$
and
$T := P^\perp$,
these~subspaces being called the algebraic and transcendental parts
of~$H^2(\frakX, \bbQ)$,
respectively.
Then~$T$
is not just a
\mbox{$\bbQ$-vector}
space,
but a pure
\mbox{$\bbQ$-Hodge}
structure \cite[Section~2.2]{De71}. The~endomorphism algebra
$\End_\Hg(T)$
is generically
just~$\bbQ$,
but may as well be a totally real field
$E \supsetneqq \bbQ$
or a CM-field~\cite[Theorems~1.6.a) and~1.5.1]{Za}.

For~an analysis of real and complex multiplications from the analytic point of view, we refer the reader to~\cite{vG}, for an interesting application of real multiplication to~\cite{Ch14}. On~the other hand, at least as far as real multiplication is concerned, the only explicit example surfaces known seem to be the ones presented by the~authors in~\cite{EJ14}. In~that article, a few
\mbox{$1$-parameter}
families were given, as well as some isolated examples, which conjecturally have~RM. For~one of these families, RM
by~$\bbQ(\sqrt{2})$
was proven, at least for countably many of its~members.


\begin{term}
\begin{iii}
\item
If~$\bbQ \subsetneqq E = \End_\Hg(T)$
is a totally real field then we say that
$\frakX$
has {\em real multiplication (RM)\/}
by~$E$.
If~$E$
is a CM-field then we say that
$\frakX$
has {\em complex multiplication (CM)\/}
by~$E$.
\item
In~either case, we call
$E = \End_\Hg(T)$
the {\em endomorphism field\/}
of~$\frakX$
(and~$T$).
\item
If~it happens that
$\bbQ \subsetneqq E' \subseteq E = \End_\Hg(T)$
then
$\frakX$
(and~$T$)
are said to be {\em acted upon\/}
by~$E'$.
\end{iii}
\end{term}

\begin{mres}
The following criterion is the main result of this article. It is a relative version of~\cite[Lemma~6.1]{EJ14}.
\end{mres}

\begin{theo}[Sufficient criterion for RM or CM in families]
\label{suff_crit}
Consider~a proper and smooth morphism\/
$\smash{\underline{q}\colon \underline{X} \to \underline{B}}$
of irreducible schemes of finite type
over\/~$\bbZ[\frac1l]$,
every fibre of which is a
$K3$~surface.
Suppose~that\/
$B$
is a normal scheme and has a\/
\mbox{$\bbQ$-rational}~point.

\begin{abc}
\item
Assume~that there exists a number
field\/~$K$
that is Galois
over\/~$\bbQ$
and a conjugacy class\/
$c$
of elements
in\/~$\Gal(K/\bbQ)$
with the property~below: For~every\/
prime
number\/~$p$
such that\/
$\Frob_p \in c$
and every\/
\mbox{$\bbF_{\!p}$-rational}
point\/
$\smash{\tau \in \underline{B}(\bbF_{\!p})}$,
the special fibre\/
$\underline{X}_\tau$
has point~count
\begin{equation}
\label{pointcount}
\#\underline{X}_\tau(\bbF_{\!p}) \equiv 1 \pmodulo p \,.
\end{equation}
Then the generic fibre
$X_\eta$
has real or complex~multiplication.
\item
Assume~that the generic fibre
$X_\eta$
has real or complex multiplication by some endomorphism
field\/~$E$.
Then,~for every complex point\/
$z \in B(\bbC)$,
the\/
$K3$~surface\/
$X_z(\bbC)$
is acted upon
by\/~$E$.
\end{abc}
\end{theo}

\begin{rems}
\begin{iii}
\item
(Cyclotomic case)
Put~$K := \bbQ(\zeta_D)$,
for an arbitrary positive
integer~$D$.
Then~$\Gal(K/\bbQ) \cong (\bbZ/D\bbZ)^*$,
which is an abelian~group. Hence,~the conjugacy classes of elements are~singletons. Moreover,~the condition
$\Frob_p \in c$~simply
means that
$p \equiv a \pmod D$,
for a certain
integer~$a$
prime
to~$D$.
\item
Even~under a congruence condition of the type described in~i), the strongest and most unlikely looking of the assumptions above is, of course, formula~(\ref{pointcount}) concerning the point~count. Note, however, that real multiplication tends to cause exactly such a behaviour~\cite[Corollary~4.13.i)]{EJ14}. In Examples~\ref{Qw2_family} and~\ref{Qw5_family}, we present explicit families of
$K3$~surfaces,
for which (\ref{pointcount}) is established by elementary~arguments.
\end{iii}
\end{rems}

\begin{rem}
In order to apply the criterion above, one further needs methods to determine the endomorphism
field~$E$
of a particular
$K3$~surface,
under the assumption that
$E \supsetneqq \bbQ$
is already~known. For~the convenience of the reader, in Section~\ref{sec_exfam_Qw5}, we add some information on how to handle the case of a quadratic~field. I.e., how to prove that
$[E:\bbQ] = 2$
and how to determine which quadratic field exactly~occurs. The method described has essentially been known before, cf.~\cite{EJ14}.
\end{rem}

\begin{rem}
We~prove Theorem~\ref{suff_crit}.a) in Section~\ref{sec_crit_fam}, as Theorem~\ref{RMCM_unspecific}, and Theorem~\ref{suff_crit}.b) in Section~\ref{sec_spec}. In~fact, for part~b), the assumptions may be somewhat weakened. Cf.~Corollary~\ref{complex} for an exact formulation.
\end{rem}

\begin{link}
Theorem~\ref{suff_crit}.a) is arithmetic in nature, taking as its main assumption condition~(\ref{pointcount}) on the numbers of points on the reductions
modulo~$p$,
for infinitely many prime
numbers~$p$.
The~link between RM and CM, i.e.\ Hodge structures, and arithmetic works as~follows.

Let~$X$
be a
$K3$~surface
over a field
$k \subset \bbC$
that is finitely generated
over~$\bbQ$.
Under~the canonical isomorphism
$\iota\colon H^2_\et(X_{\overline{k}}, \bbQ_l(1)) \to H^2(X(\bbC), \bbQ) \!\otimes_\bbQ\! \bbQ_l(1)$
between \'etale and complex cohomology \cite[Expos\'e~11, Th\'eor\`eme~4.4.iii)]{SGA4}, the algebraic classes and the cup product pairing are~respected. Thus,
$T  \!\otimes_\bbQ\! \bbQ_l(1)$,
for
$T \subset H^2(Z(\bbC), \bbQ)$
the transcendental part, gets identified with
$\calT \subset H^2_\et(X_{\overline{k}}, \bbQ_l(1))$,
the transcendental part of \'etale~cohomology. As~a consequence of this, by transport of structure, one has an~operation
\begin{equation}
\label{op_on_etale}
E = \End_\Hg(T) \hookrightarrow \End(\calT)
\end{equation}
of~$E = \End_\Hg(T)$
on~$\calT$.

It~is well-known \cite{Za,Ta90,Ta95,An} that the neutral component of the algebraic monodromy group
of~$\calT$
is given by
$\smash{\MG_{\calT, k, l}^0 = (C_{\GO(\calT)}(E))^0}$.
Cf.\ Example~\ref{alg_mon_grp}.v). I.e., except for the case of geometric Picard
rank~$20$,
in which
$\GO^0(\calT)$
is abelian and CM is automatic \cite[Remark~3.3.10]{Hu}, one has
$\smash{\MG_{\calT, k, l}^0 \subsetneqq \GO^0(\calT)}$
if and only
if~$\smash{E \supsetneqq \bbQ}$.
\end{link}

\begin{fres}
We~reverse this result in Theorem~\ref{Zarhin_rev} and~prove
$$C_{\End(\calT)}(\MG_{\calT, k, l}^0) = E \!\otimes_\bbQ\! \bbQ_l \,.$$
%
Thus,~the endomorphism field
$E$
is determined by the algebraic monodromy group, at least up to arithmetic equivalence~\cite{Pe}. This fixes the degree
of~$E$
and in many situations
$E$
itself, for example when
$E$
is Galois over
$\bbQ$
or when
$[E\!:\!\bbQ] < 7$~\cite{BdS}. 
It~turns out (cf.\ Corollary~\ref{E_indep}) that there is another situation, in which
$E$~is
independent
of~$k \hookrightarrow \bbC$,
namely when the base field
$k$
is primary
over~$\bbQ$,
i.e.\ does not contain any proper algebraic extension
of~$\bbQ$.
\end{fres}

\begin{mterm}
Let~$k$
be a field that is finitely generated
over~$\bbQ$
and
$X$
a
$K3$
sur\-face
over~$k$.
Assume~that
$k$
is primary
over~$\bbQ$
or
that~$\rk\Pic X_{\overline{k}} > 15$.

\begin{iii}
\item
We~say that
$X$~has
real or complex multiplication if
$X(\bbC)$~has.
(This terminology is used in Theorem~\ref{suff_crit}.a) for the generic fibre
$X_\eta$.)
\item
Similarly,~we shall fell free to speak of the endomorphism field of
$X$,
instead
of~$X(\bbC)$.
\end{iii}
\end{mterm}

\begin{appl}
As~an application, in Section~\ref{sec_exfam_Qw2}, we return to the family from~\cite{EJ14} and prove that actually every member is acted upon
by~$\smash{\bbQ(\sqrt{2})}$.
Moreover,~we present a new family, all members of which are acted upon
by~$\smash{\bbQ(\sqrt{5})}$.
\end{appl}

\begin{ex}[An explicit family of
$K3$~surfaces
with RM
by~$\bbQ(\sqrt{5})$]
\label{Qw5_family}
Let
$\smash{q\colon X \to B}$,
for
$\smash{B := \Spec \bbQ[t,\!\frac1{(t-1)(t^4-t^3+t^2-t+1)}] \subset \Ab^1_\bbQ}$,
be the family of
$K3$~surfaces
that is fibre-by-fibre the minimal desingularisation of the double cover
of~$\Pb^2$,
given~by
\begin{eqnarray}
\label{Qw5}
w^2 &=&  y (x - 2(t\!-\!1)y - tz) \\[-1mm]
&& (x^4 + x^3y - x^3z + x^2y^2 - 2x^2yz + x^2z^2 + xy^3 - 3xy^2z - 2xyz^2 - xz^3 + y^4 \nonumber \\[-1mm]
&& \hspace{9.2cm} {}+ y^3z + y^2z^2 + yz^3 + z^4) \, . \nonumber
\end{eqnarray}
\begin{iii}
\item
Then the generic fibre
$X_\eta$
of~$q$
is of geometric Picard
rank~$16$.
\item
The endomorphism field of
$X_\eta$
is~$\smash{\bbQ(\sqrt{5})}$.
\item
For~every
$\theta \in B(\bbC)$,
the transcendental part
$\smash{T \subset H^2(X_\theta(\bbC), \bbQ)}$
of the cohomology of the fibre
$X_\theta$
is acted upon
by~$\smash{\bbQ(\sqrt{5})}$.
\item
Let the complex point
$\theta \in B(\bbC)$
be of the kind that the fibre
$X_\theta$
has Picard
rank~$16$.
Then
$X_\theta$
has real multiplication
by~$\smash{\bbQ(\sqrt{5})}$.
\end{iii}\smallskip

\noindent
{\bf Proof.}
We~prove these results in Section~\ref{sec_exfam_Qw5}.
\eop
\end{ex}

\begin{tate}
Our arguments make use of the Tate conjecture, in situations where it is known to be true. More concretely, what we use is the~following.
\end{tate}

\begin{facs}[Known cases of the Tate conjecture]
\label{Tate}
Let\/~$k$
be field that~is

\begin{abc}
\item
finitely generated
over\/~$\bbQ$
or
\item
a finite field
\end{abc}
and\/
$X$
a\/~$K3$~surface
over\/~$k$.
Then~the subspace\/
$\smash{H^2_\et(X_{\overline{k}}, \bbQ_l(1))^{\Gal(\overline{k}/k)}}$
of invariants coincides with\/
$c_1(\Pic X) \!\otimes_\bbZ\! \bbQ_l \subset H^2_\et(X_{\overline{k}}, \bbQ_l(1))$.\medskip

\noindent
{\bf Proof.}
{\em
a)
This~is due to Y.\ Andr\'e~\cite[Section~6.2]{An}. More~recent developments are described in~\cite[Proposition~9.2]{Mo}.\smallskip

\noindent
b)
In this situation, the result is known due to the combined work of several people, most notably F.~Charles \cite{Ch13}, M.~Lieb\-lich, D.\ Maulik, and A.\ Snowden~\cite{LMS}, K.~Madapusi Pera \cite{MP}, as well as W.~Kim and K.\ Madapusi Pera \cite{KM}.
}
\eop
\end{facs}

\begin{conv}
We~follow standard conventions and use standard notation from Algebra and Algebraic~Geometry. More~specifically,

\begin{iii}
\item
For a field
$k$,
we denote
by~$\overline{k}$
a separable closure. For~a point
$t\colon \Spec k \to X$,
we write
$\smash{\overline{t}\colon \Spec\overline{k} \to \Spec k \to X}$
for the resulting geometric~point.
\item
We usually denote the generic point on a connected scheme
by~$\eta$.
\item
We~say that a proper scheme
$X$
over a number field
$k$
has good reduction at a prime ideal
$\frakp$
of~$k$,
if there exists a proper model
$\underline{X}$
over the integer ring
$\calO_k \subset k$
that is smooth
above~$\frakp$.
\item
When
$\underline{B}, \underline{X}, \ldots$
is a scheme of finite type
over~$\Spec\bbZ$,
we write
$B, X$
etc.\ for its generic fibre. E.g.,
$B := \underline{B} \times_{\Spec\bbZ} \Spec\bbQ$.
\item
For~a
$K3$~surface
$X$
over a field
$k \subset \bbC$
that is finitely generated
over~$\bbQ$,
the canonical comparison isomorphism \cite[Expos\'e~11, Th\'eor\`eme~4.4.iii)]{SGA4} induces an isomorphism
$\calT \cong T \!\otimes_\bbQ\! \bbQ_l(1)$,
for
$T \subset H^2(X(\bbC), \bbQ)$
and
$\calT \subset H^2_\et(X_{\overline{k}}, \bbQ_l(1))$
the transcendental~parts. We~identify the two
\mbox{$\bbQ_l$-vector}
spaces and consider
$T$~as
a subset
of~$\calT$.
\item
For an algebraic
group~$G$,
we denote
by~$G^0$
its neutral component with respect to the Zariski~topology. Similarly for the set of all
\mbox{$\bbQ_l$-rational}
points on an algebraic group defined
over~$\bbQ_l$.
\item
We denote by
$\GO_n$
the linear algebraic group
$\bbG_m \!\cdot\! \O_n$.
In~characteristic different
from~$2$,
$\GO_n$
is irreducible for
$n$~odd
and has two components for
$n$~even.
In~any case, one has
$\smash{\GO^0_n = \bbG_m \!\cdot\! \SO_n}$.

When
$\calT$
is a finite-dimensional
\mbox{$\bbQ_l$-vector}
space equipped with a non-degenerate symmetric bilinear form, we use the notation
$\GO(\calT)$
for the group of all orthogonal similitudes
of~$\calT$.
\item
When
$G,H \subseteq A$
are groups or algebras that are contained in an
algebra~$A$
then
$C_H(G) := \{h \in H \mid \forall g \in G\colon hg=gh \}$
denotes the centraliser of
$G$
in~$H$.
\item
Except~for Corollary~\ref{E_indep_el}, we fix a prime number
$l$
throughout the article.
\end{iii}
\end{conv}

\section{Algebraic monodromy groups}

Algebraic monodromy groups are the main tool that is used in the present~article. The~purpose of this section is primarily to recall the relevant facts and to fix~notation.

Let~$B$
be an arbitrary connected scheme, on which a geometric point
$\smash{\overline{s}}$
is fixed as the base~point. As~is well-known \cite[Expos\'e~VI, Lemme~1.2.4.2]{SGA5}, associated with any
\mbox{$\bbQ_l$-sheaf}
$\calQ$
on~$B$
that is twisted-constant with respect to the \'etale topology and of finite rank, one has a continuous~representation
\begin{equation}
\label{rep_fundgr}
\varrho^\calQ_{\overline{s}}\colon \pi_1^\et(B,\overline{s}) \longrightarrow \GL(\calQ_{\overline{s}})
\end{equation}
of the \'etale fundamental group
of~$B$.

\begin{defi}
The~Zariski closure of the image of
$\varrho^\calQ_{\overline{s}}$
is called the {\em algebraic monodromy group\/}
$\MG_{\calQ, B, l}$
of~$\calQ$.
This~is the set of all
\mbox{$\bbQ_l$-rational}
points on an algebraic group defined
over~$\bbQ_l$.
The~algebraic group is possibly~disconnected.
\end{defi}

\begin{rem}
\label{indep_basepoint}
Let
$\overline{s}'$
be another geometric point
on~$B$.
One~may then choose an \'etale path 
$\smash{\gamma \in \pi_1^\et(B, \overline{s}, \overline{s}')}$,
which yields an isomorphism
$i_\gamma\colon \calQ_{\overline{s}} \to \calQ_{\overline{s}'}$.
The~diagram
$$
\xymatrixcolsep{2.7mm}
\xymatrix{
\pi_1^\et(B, \overline{s}) \ar@{->}[rr]^{\varrho^\calQ_{\overline{s}}} \ar@{->}[d]_{\sigma \mapsto \gamma\sigma\gamma^{-1}} && \GL(\calQ_{\overline{s}}) \ar@{->}[d]^{M \mapsto i_\gamma M i_\gamma^{-1}} \\
\pi_1^\et(B, \overline{s}') \ar@{->}[rr]^{\varrho^\calQ_{\overline{s}'}} && \GL(\calQ_{\overline{s}'})
}
$$
then~commutes. I.e.,
$\gamma$
induces an isomorphism between
$\GL(\calQ_{\overline{s}})$
and
$\GL(\calQ_{\overline{s}'})$
that maps the two instances of
$\MG_{\calQ, B, l}$
onto each other.
\end{rem}

\begin{exs}
\label{alg_mon_grp}
\begin{iii}
\item
Let
$B = \Spec k$,
for
$k$
a~field. Then~(\ref{rep_fundgr}) specialises to a representation of
$\smash{\Gal(\overline{k}/k) = \pi_1^\et(\Spec k,\Spec\overline{k})}$.
The algebraic monodromy group
$\MG_{\calQ, k, l}$
of~$\calQ$
is the Zariski closure of the image of
$\smash{\Gal(\overline{k}/k)}$
in
$\GL(\calQ_{\overline{k}})$.
\item
When~$B$
is an arbitrary scheme of residue characteristics different
from~$l$
and
$q\colon X \to B$
a smooth and proper morphism then the proper and smooth base change theorems \cite[Expos\'e~XVI, Corollaire 2.2]{SGA4} imply that the higher direct image sheaves
$\smash{R^i q_* \bbQ_l(j)}$
are twisted-constant
\mbox{$\bbQ_l$-sheaves},
for all
$i, j \in \bbN$.
One has the algebraic monodromy group
$\smash{\MG_{R^i q_* \bbQ_l(j), B, l} \subseteq \GL((R^i q_* \bbQ_l(j))_{\overline{s}}) = \GL(H^i_\et(X_{\overline{s}}, \bbQ_l(j)))}$.
\item
Suppose that the fibres
of~$q$
are of
dimension~$i$,
for an even
integer~$i$. Then~the stalks
of~$\smash{R^i q_* \bbQ_l(i/2)}$
at the geometric point
$\overline{s}$
are equipped with the symmetric pairing, induced by cup product and Poincar\'e duality. As~a consequence of~this,
$$\MG_{R^i q_* \bbQ_l(i/2), B, l} \subseteq \GO((R^i q_* \bbQ_l(i/2))_{\overline{s}}) = \GO(H^i_\et(X_{\overline{s}}, \bbQ_l(i/2))) \, .$$
\item
Assume, in addition, that there is given a decomposition
$\smash{R^i q_* \bbQ_l(i/2) = \calP \oplus \calT}$
into two twisted-constant subsheaves, whereas the restriction of the cup product pairing to
$\calP_{\overline{s}}$
is non-degenerate. Then the same is true for
$\calT_{\overline{s}}$
and the algebraic monodromy group
$\smash{\MG_{\calT, B, l}}$
is contained in
$\smash{\GO(\calT_{\overline{s}})}$.
\item
(Zarhin, Tankeev, Andr\'e)
Consider the situation that
$B = \Spec k$,
for
$k$
a field that is finitely generated
over~$\bbQ$,
$X$
is a
$K3$~surface
over~$k$,
and
$\calT \subset H^2_\et(X_{\overline{k}}, \bbQ_l(1))$
is the transcendental~part.

Then~the neutral component of the algebraic monodromy group
of~$\calT$
is equal to the neutral component of the centraliser
of~$E$
in~$\GO(\calT)$,
\begin{equation}
\label{monodromygrp_centr}
\MG_{\calT, k, l}^0 = (C_{\GO(\calT)}(E))^0 \,.
\end{equation}
Here,~the endomorphism field
$E = \End_\Hg(T)$
is considered as being contained in
$\End(\calT)$
via the operation~(\ref{op_on_etale}).

Indeed,~(\ref{monodromygrp_centr}) follows from Yu.\,G.\ Zarhin's explicit description of the Mumford-Tate group in the case of a complex
$K3$~surface
\cite[Theorem~2.2.1]{Za}, together with the Mum\-ford-Tate~conjecture, cf.~\cite[Theorem~13]{Ch14}. The~Mumford-Tate~conjecture was proven for
$K3$~surfaces
over number fields by S.\,G.\ Tankeev \cite{Ta90,Ta95} and over arbitrary finitely generated extensions
of~$\bbQ$
by Y.\ Andr\'e \cite[Th\'eor\`eme~8.2]{An}. Cf.~\cite[Theorem~1.1]{Com} for recent developments.
\end{iii}
\end{exs}

\begin{ttt}[Base change]
\label{bc_spec}
Let
$i\colon B' \to B$
be a morphism of connected schemes and
$\overline{s}$
any geometric point
on~$B'$.
Then~$i$
induces an isomorphism
$\smash{(i^* \calQ)_{\overline{s}} \cong \calQ_{\overline{i(s)}}}$
and a homomorphism
$\smash{i_\#\colon \pi_1^\et(B', \overline{s}) \to \pi_1^\et(B, \overline{i(s)})}$,
via which one has a natural~inclusion
$$\MG_{i^* \calQ, B', l} \hookrightarrow \MG_{\calQ, B, l}$$
of subgroups of
$\smash{\GL((i^* \calQ)_{\overline{s}}) \cong \GL(\calQ_{\overline{i(s)}})}$.

\begin{iii}
\item
This applies, of course, when
$s\colon \Spec k \to B$
is a point,
$$\MG_{\calQ |_s, k, l} \hookrightarrow \MG_{\calQ, B, l} \,.$$
Note~here that, in view of Remark~\ref{indep_basepoint}, one may assume the base point
on~$B$
to be chosen as an extension
of~$s$.
\item
In the particular case that
$B$
is normal and locally Noetherian and that
$\eta \in B$
is the generic point, the natural inclusion is actually a~bijection,
$$\MG_{\calQ |_\eta, k(\eta), l} = \MG_{\calQ, B, l} \,.$$
Indeed,~the homomorphism
$i_\#$
is then surjective, according to~\cite[Expos\'e~V, Proposition~8.2]{SGA1}.
\end{iii}
\end{ttt}

\section{A criterion for RM or CM in families}
\label{sec_crit_fam}

Let~$\underline{B}$
be a connected scheme of finite type
over~$\bbZ$,
$q$
an arbitrary prime power, and
$\tau\colon \Spec\bbF_{\!q} \to \underline{B}$
an arbitrary closed point, defined
over~$\bbF_{\!q}$.
Then,~by functoriality, there is the natural~homomorphism
$$\tau_\#\colon \Gal(\overline\bbF_{\!q}/\bbF_{\!q}) = \pi_1^\et(\Spec\bbF_{\!q}, \Spec\overline\bbF_{\!q}) \to \pi_1^\et(\underline{B}, \overline\tau) \,,$$
so
$\tau$
defines a unique
element~$\smash{\F_{\!\tau} \in \pi_1^\et(\underline{B}, \overline\tau)}$
being the image of the canonical generator
$\smash{\Frob \in \Gal(\overline\bbF_{\!q}/\bbF_{\!q})}$.
As
$\smash{\underline{B}}$
is connected, an \'etale path
$\smash{\gamma \in \pi_1^\et(\underline{B}, \overline\tau, \overline\eta)}$
may be~chosen. This~yields an isomorphism
$$\pi_1^\et(\underline{B}, \overline\tau) \longrightarrow \pi_1^\et(\underline{B}, \overline\eta), \quad \sigma \mapsto \gamma\sigma\gamma^{-1} \,,$$
under which
$\F_{\!\tau}$
is sent to the {\em Frobenius
element\/}~$\smash{\Frob_\tau \in \pi_1^\et(\underline{B}, \overline\eta)}$.
The~Frobenius element
$\Frob_\tau$
is unique only up to conjugation, as one may choose various \'etale paths.

\begin{rem}
In the particular case that
$\smash{\underline{B} := \Spec\bbZ[\frac1D]}$,
for some integer
$D \neq 0$,
one has the Frobenius element
$\smash{\Frob_p \in \pi_1^\et(\Spec\bbZ[\frac1D], \overline\eta)}$,
for any prime
number~$p \,\notd D$.
It~causes an automorphism of every Galois \'etale covering
of~$\smash{\Spec\bbZ[\frac1D]}$
and, in particular, an element
in~$\Gal(K/\bbQ)$,
as long as
$K$
is a number field that is normal
over~$\bbQ$
and ramified only at primes
dividing~$D$.
This~is the classical Frobenius element from Algebraic Number~Theory.
\end{rem}

Our argumentation in this section essentially relies on the Cebotarev density theorem in the extended version, due to \mbox{J.-P.}\ Serre~\cite{Serr}. This~is the following~result.

\begin{prop}[J.-P.\ Serre]
\label{Cebotarev_Serre}
Let\/~$\underline{B}$
be an irreducible scheme of finite type
over\/~$\bbZ$
and\/
$U \subset \pi_1^\et(\underline{B}, \overline\eta)$
a normal subgroup of finite~index.

\begin{abc}
\item
Then~the Frobenius elements\/
$\Frob_\tau \in \pi_1^\et(\underline{B}, \overline\eta)$
of the closed points\/
$\tau \in \underline{B}$
take every conjugacy class of elements
of\/~$\pi_1^\et(\underline{B}, .)/U$.
\item
Suppose, in addition, that\/
$\underline{B}$
is flat
over\/~$\bbZ$.
Then~the Frobenius elements\/
$\Frob_\tau$
of the closed points\/
$\tau \in \underline{B}$
defined over a prime field already take every conjugacy class of elements
of\/~$\pi_1^\et(\underline{B}, .)/U$.
\end{abc}\smallskip

\noindent
{\bf Proof.}
{\em
a)
directly follows from \cite[Theorem~7]{Serr}. In~fact, J.-P.\ Serre shows for every conjugacy class that the suitable closed points are of positive Dirichlet~density.\smallskip

\noindent
b)
Thus,~to establish~b), it suffices to show that the closed points defined over non-prime fields form a set of Dirichlet density~zero. For~this, let us write
$n := \dim \underline{B}$.
Then~the special fibre
over~$\bbF_{\!p}$
is of dimension
$(n-1)$.
Hence,~according to the Lang--Weil estimates~\cite[Theorem~1]{LW}, there is some constant
$C\in\bbR$
such that
$\#\underline{B}(\bbF_{\!q}) \leq Cq^{n-1}$,
for every prime
power~$q$.
Therefore,
$$\sum_{\atop{p\,\text{prime}}{k\geq2, \tau\in\underline{B}(\bbF_{\!p^k})}} \!\!\!\!\frac1{N(\tau)^n} \leq C \sum_{\atop{p\,\text{prime}}{k\geq2}} \frac{(p^k)^{n-1}}{(p^k)^n} = C \sum_{\atop{p\,\text{prime}}{k\geq2}} \frac1{p^k} < \infty \,,$$
which implies the claim, cf.~\cite[formula~(19)]{Serr}.
}
\eop
\end{prop}

\begin{theo}[Strict inclusion for the algebraic monodromy group]
\label{pc_yields_monodr}
Let\/~$D \neq 0$
be an integer and\/
$\smash{\underline{q}\colon \underline{X} \to \underline{B}}$
a proper and smooth morphism of irreducible schemes of finite type
over\/~$\bbZ[\frac1{lD}]$,
every fibre of which is a\/
$K3$~surface.
Suppose~that\/
$B$~has
a\/
\mbox{$\bbQ$-rational}
point.\smallskip

\noindent
Let, moreover,
$\smash{\calT \subset R^2 q_* \bbQ_l(1)}$
be a twisted-constant sheaf of
rank~$\geq \!3$
of the kind that

\begin{iii}
\item[$\bullet$ ]
there is a decomposition\/
$\smash{R^2 q_* \bbQ_l(1) = (\calP \!\otimes_\bbZ\! \bbQ_l) \oplus \calT}$,
for\/
$\calP$
a locally constant\/
\mbox{$\bbZ$-sheaf,}
\item[$\bullet$ ]
the restriction of the cup product pairing to\/
$\calT_{\overline{\eta}}$
is non-degenerate.
\end{iii}

\noindent
Finally,~assume that there exists a number
field\/~$K$
of discriminant\/~$D$
that is Galois
over\/~$\bbQ$
and a conjugacy class\/
$c$
of elements
in\/~$\Gal(K/\bbQ)$
with the property~below: For~every~prime
number\/~$p$
such that\/
$\Frob_p \in c$
and every\/
\mbox{$\bbF_{\!p}$-rational}
point\/
$\smash{\tau \in \underline{B}(\bbF_{\!p})}$,
the special fibre\/
$\underline{X}_\tau$
has point~count
\begin{equation}
\label{pc}
\#\underline{X}_\tau(\bbF_{\!p}) \equiv 1 \pmodulo p \,.
\end{equation}
Then~the strict~inclusion\/
$\MG^0_{\calT, \underline{B}, l} \subsetneqq \GO^0(\calT_{\overline{\eta}})$~holds.
\end{theo}

\begin{coro}[Algebraic monodromy group of the generic fibre]%
\label{amg_fibres}
Let\/
$\smash{\underline{q}\colon \underline{X} \to \underline{B}}$
and\/
$\smash{\calT \subset R^2 q_* \bbQ_l(1)}$
be as in Theorem~\ref{pc_yields_monodr}.
Then
$$\smash{\MG^0_{\calT_{\overline\eta, k(\eta), l}} \subsetneqq \GO^0(\calT_{\overline\eta}) \,.}\smallskip$$

\noindent
{\bf Proof.}
{\em
This follows directly from Theorem~\ref{pc_yields_monodr}, together with \ref{bc_spec}.i).
}
\eop
\end{coro}

\begin{theo}[Sufficient criterion for the generic fibre to have RM or CM with an unspecified endomorphism field]
\label{RMCM_unspecific}
Let\/
$\smash{\underline{q}\colon \underline{X} \to \underline{B}}$
be a proper and smooth morphism of irreducible schemes of finite type
over\/~$\bbZ[\frac1l]$,
every fibre of which is a\/
$K3$~surface.
Suppose~that
$B$~is
a normal scheme and has
a\/
\mbox{$\bbQ$-rational}
point.\smallskip

\noindent
Assume~that there exists a number
field\/~$K$
that is Galois
over\/~$\bbQ$
and a conjugacy class\/
$c$
of elements
in\/~$\Gal(K/\bbQ)$
with the property~below: For~every~prime
number\/~$p$
such that\/
$\Frob_p \in c$
and every\/
\mbox{$\bbF_{\!p}$-rational}
point\/
$\smash{\tau \in \underline{B}(\bbF_{\!p})}$,
the special fibre\/
$\underline{X}_\tau$
has point~count
$$\#\underline{X}_\tau(\bbF_{\!p}) \equiv 1 \pmodulo p \,.$$
Then~the\/
$K3$~surface\/
$X_\eta(\bbC) = (X_\eta \times_{\Spec k(\eta)} \Spec\bbC)(\bbC)$
has real or complex multiplication,~for every embedding\/
$k(\eta) \hookrightarrow \bbC$.\medskip

\noindent
{\bf Proof.}
{\em
Without restriction,
$\smash{\underline{q}}$
is a morphism of
\mbox{$\smash{\bbZ[\frac1{lD}]}$-schemes},
for~$D$
the discriminant of the field
extension~$K/\bbQ$.
Indeed,
$$\textstyle \underline{q} \!\times_{\Spec\bbZ[\frac1l]}\! \Spec\bbZ[\frac1{lD}]\colon \underline{X} \!\times_{\Spec\bbZ[\frac1l]}\! \Spec\bbZ[\frac1{lD}] \longrightarrow \underline{B} \!\times_{\Spec\bbZ[\frac1l]}\! \Spec\bbZ[\frac1{lD}]$$
still fulfils all the assumptions~made.

Put~$\smash{\calR := R^2 \underline{q}_* \bbQ_l(1)}$
and
$\smash{\calP_{\overline\eta} := c_1(\Pic X_{\overline\eta}) \subset H^2_\et(X_{\overline\eta}, \bbQ_l(1)) = \calR_{\overline\eta}}$.
Then~$\calR$
is a twisted-constant
\mbox{$\bbQ_l$-sheaf}
on~$\smash{\underline{B}}$,
due to the smooth and proper base change theorems \cite[Expos\'e~XVI, Corollaire~2.2]{SGA4}.
Moreover,~$\smash{\calP_{\overline\eta}}$
is clearly
$\smash{\Gal(\overline{k(\eta)}/k(\eta))}$-invariant
and stabilised by an open subgroup of finite index, as every invertible sheaf is defined over a finite extension
of~$k(\eta)$.
Therefore,~Lemma~\ref{sheaf_ext} applies and shows that
$\smash{\calP_{\overline\eta}}$
extends to a locally constant
\mbox{$\bbZ$-sheaf}~$\calP$
on~$\smash{\underline{B}}$.

Finally,~write
$\calT := (\calP \!\otimes_\bbZ\! \bbQ_l)^\perp$.
If~$\rk\calT \geq 3$
then all assumptions of Theorem~\ref{pc_yields_monodr} are~satisfied. Indeed, the cup product pairing
on~$\calP$
is non-degenerate and this implies the same
for~$\calT$.
In view of formula~(\ref{monodromygrp_centr}), the assertion is then a direct consequence of Corollary~\ref{amg_fibres}.
Otherwise,~one has
$\rk\calT = 2$
and hence
$X_\eta$
is of geometric Picard
rank~$20$.
In~this case,
$X_\eta(\bbC)$
is known to have CM \cite[Remark~3.3.10]{Hu}.
}
\eop
\end{theo}

\begin{lem}
\label{sheaf_ext}
Let\/~$A$
be a normal scheme that is connected and locally Noetherian and
let\/~$\calR$
be a twisted-constant\/
\mbox{$\bbQ_l$-sheaf}
on\/~$A$.
Moreover,~let\/
$S \hookrightarrow S \!\otimes_\bbZ\! \bbQ_l \subseteq \calR_{\overline\eta}$
be a free\/
\mbox{$\bbZ$-module}
of finite rank that is invariant under the\/
\mbox{$\smash{\Gal(\overline{k(\eta)}/k(\eta))}$-operation}
and stabilised by an open subgroup of finite index
in\/~$\smash{\Gal(\overline{k(\eta)}/k(\eta))}$.\smallskip

\noindent
Then\/~$S$
extends to the whole
of\/~$A$
as a locally constant\/
\mbox{$\bbZ$-sheaf\/}~$\calS$.
Moreover, one has\/
$\calS \!\otimes_\bbZ\! \bbQ_l \subseteq \calR$.\medskip

\noindent
{\bf Proof.}
{\em
The~assumptions made
on~$A$
are enough to imply that the \'etale fundamental group
$\pi_1^\et(A,.)$
is a quotient of
$\smash{\Gal(\overline{k(\eta)}/k(\eta))}$
\cite[Expos\'e~V, Proposition~8.2]{SGA1}.
Write~$\smash{\pi_1^\et(A,.) = \Gal(\overline{k(\eta)}/k(\eta))/H}$.
The~equivalence of categories between twisted-constant
\mbox{$\bbQ_l$-sheaves}
on~$A$
and
\mbox{$\bbQ_l$-vector}
spaces being continuously acted upon
by~$\pi_1^\et(A,.)$
\cite[Expos\'e~VI, Lemme~1.2.4.2]{SGA5} therefore shows that
$\calR_{\overline\eta}$
is actually a
\mbox{$\smash{\Gal(\overline{k(\eta)}/k(\eta))/H}$-module}.
In~particular,
$H$
operates trivially
on~$S$.

Thus,
$S$
is acted upon
by~$\smash{\Gal(\overline{k(\eta)}/k(\eta))/H = \pi_1^\et(A,.)}$.
The~assumption implies that
$S$
is, furthermore, stabilised by an open subgroup
$\smash{K \subset \pi_1^\et(A,.)}$
of finite~index. Consequently,~there exists an \'etale covering
$\smash{\widetilde{A} \to A}$,
which we may assume to be Galois, on which
$S$
defines a constant
\mbox{$\bbZ$-sheaf}
$\smash{\widetilde\calS}$,
together with an operation of the finite group
$\smash{\Gal(\widetilde{A}/A) = \pi_1^\et(A,.)/K}$.

Since every orbit of a finite group is finite, 
$\smash{\widetilde\calS}$~is
represented \cite[Expos\'e~VII, Section~2.a)]{SGA4} by an infinite disjoint
union~$\frakA$
of trivial, finite \'etale coverings, i.e.\ such of type
$\smash{\widetilde{A} \sqcup\cdots\sqcup \widetilde{A}}$,
each of which is acted upon
by~$\pi_1^\et(A,.)/K$.
According~to \cite[Expos\'e~IX, Proposition~4.1]{SGA1},
$\frakA$~descends
to an infinite disjoint union of finite \'etale coverings
of~$A$.
This~union, finally, represents a
sheaf~$\calS$
of sets
on~$A$.

The
sheaf~$\calS$
is locally constant, since it is trivialised by the \'etale covering
$\smash{\widetilde{A} \to A}$.
As~the group structure descends, too, the proof is~complete.
}
\eop
\end{lem}

\begin{proo1}
{\em First step.}
For~an arbitrary
$\bbF_{\!p}$-rational
point
$\tau\in\underline{B}(\bbF_{\!p})$,
for
$p$
of the kind that
$\Frob_p \in c$,
one has
$\smash{\Tr(\varrho^\calT_{\overline\eta}(\Frob_\tau)) \in [-22,22] \cap \bbZ}$.

\noindent
The Lefschetz trace formula \cite[Expos\'e~XII, 6.3 and Exemple 7.3]{SGA5} yields that
\begin{align*}
\#\underline{X}_\tau(\bbF_{\!p}) &= p^2 + \Tr(\Frob\colon H^2_\et(X_{\overline\tau}, \bbQ_l(1)) \to H^2_\et(X_{\overline\tau}, \bbQ_l(1))) p + 1 \\
 &= p^2 + \Tr(\F_{\!\tau}\colon (R^2 q_* \bbQ_l(1))_{\overline\tau} \to (R^2 q_* \bbQ_l(1))_{\overline\tau}) p + 1 \,.
\end{align*}
Note~the factor
$p$
that is a result of the Tate~twist. Consequently,~Assumption~(\ref{pc}) means nothing~but
$$\Tr \!\big(\! \F_{\!\tau}\colon (R^2 q_* \bbQ_l(1))_{\overline\tau} \to (R^2 q_* \bbQ_l(1))_{\overline\tau} \big) \in\bbZ \,.$$
But, via any \'etale
path~$\smash{\gamma \in \pi_1^\et(\underline{B}, \overline\tau, \overline\eta)}$,
\begin{align*}
\Tr(\varrho^{R^2 q_* \bbQ_l(1)}_{\overline\eta}(\Frob_\tau)) &= \Tr(\Frob_\tau\colon (R^2 q_* \bbQ_l(1))_{\overline\eta} \to (R^2 q_* \bbQ_l(1))_{\overline\eta}) \\
 &= \Tr(\F_{\!\tau}\colon (R^2 q_* \bbQ_l(1))_{\overline\tau} \to (R^2 q_* \bbQ_l(1))_{\overline\tau}) \,.
\end{align*}
On the other hand,
$\smash{\Tr(\varrho^{\calP \otimes_\bbZ \bbQ_l}_{\overline\eta}(\Frob_\tau)) = \Tr(\varrho^\calP_{\overline\eta}(\Frob_\tau)) \in \bbZ}$.
Note~here that
$\Frob_\tau$
operates already on the
\mbox{$\bbZ$-sheaf}~$\calP$.
Consequently,~one has
$\smash{\Tr(\varrho^\calT_{\overline\eta}(\Frob_\tau)) \in \bbZ}$
too.

The~argument above also shows that
$\smash{\Tr(\varrho^\calT_{\overline\eta}(\Frob_\tau))}$
is the same as the sum of the eigenvalues of
$\Frob$,
operating on
$\smash{\calT_{\overline\tau} \subset (R^2 q_* \bbQ_l(1))_{\overline\tau} = H^2_\et(X_{\overline\tau}, \bbQ_l(1))}$.
Since these are all of absolute
value~$1$,
according to the Weil conjectures, proven by P.\ Deligne \cite[Th\'eo\-r\`eme 8.1]{De74}, and the vector space to the right is of
dimension~$22$,
the claim follows.\smallskip

\noindent
{\em Second step.}
A finite-index subgroup
of~$\smash{\pi_1^\et(\underline{B}, \overline\eta)}$.

\noindent
The~homomorphism
$\smash{\pr_\#\colon \pi_1^\et(\underline{B}, \overline\eta) \to \pi_1^\et(\Spec\bbZ[\frac1{lD}], \Spec\overline\bbQ)}$,
induced by the structural morphism
$\smash{\pr\colon \underline B \to \Spec\bbZ[\frac1{lD}]}$,
is~surjective. Indeed,~for
$t\colon \Spec\bbQ \to \underline{B}$
any
\mbox{$\bbQ$-rational}
point, it is sufficient to show that
$\smash{\pr_\# \!\circ\! t_\# = (\pr \!\circ\! t)_\#}$
is~surjective.
But~as
\mbox{$\smash{\pr \!\circ\! t\colon \Spec\bbQ \to \Spec\bbZ[\frac1{lD}]}$}
is just the embedding of the generic point, this is true, due to \cite[Expos\'e~V, Proposition~8.2]{SGA1}.

Furthermore, there is the homomorphism
$$\textstyle \chi\colon \pi_1^\et(\Spec\bbZ[\frac1{lD}], \Spec\overline\bbQ) \to \Gal(K/\bbQ) \,,$$
induced by the operation on the
\mbox{$[K\!:\!\bbQ]$-fold}
\'etale covering
$$\textstyle W := \Spec\calO_K[\frac1{lD}] \to \Spec\bbZ[\frac1{lD}] \,.$$
The
homomorphism~$\chi$
is epic, too, as
$W$
is a connected scheme.
Put~$\chi_K := \chi \!\circ\! \pr_\#$,
$$U := \ker \chi_K \quad \text{and}\quad U^{(c)} := \chi_K^{-1}(c) \,.$$
Then~$\smash{U \subset \pi_1^\et(\underline{B}, \overline\eta)}$
is an open subgroup of
index~$[K\!:\!\bbQ]$
and
$\smash{U^{(c)} \subset \pi_1^\et(\underline{B}, \overline\eta)}$
is a non-empty union of cosets
modulo~$U$.\medskip

\noindent
{\em Third step.}
For every\/
$\sigma \in U^{(c)}$,
one has
$\smash{\Tr(\varrho^\calT_{\overline\eta}(\sigma)) \in [-22,22] \cap \bbZ}$.

\noindent
Assume~the contrary. Then
$\smash{\Tr(\varrho^\calT_{\overline\eta}(\sigma)) \not\in \{-22, \ldots, 22\} \subset \bbQ_l}$,
for one particular such
element~$\sigma$.

In~order to analyse this assumption further, let us consider a torsion-free
\mbox{$\bbZ_l$-sheaf}
$\calT^{(\bbZ_l)}$
underlying~$\calT$.
Such~a sheaf must exist for very general reasons \cite[Expos\'e~6, D\'efinition~1.4.3]{SGA5}. Consequently, there is a continuous representation
$$\varrho^{\calT^{(\bbZ_l)}}_{\overline\eta}\colon \pi_1^\et(\underline{B}, \overline\eta) \to \GL(\calT^{(\bbZ_l)}_{\overline\eta})$$
underlying~$\varrho^\calT_{\overline\eta}$.
In~particular, we have\vspace{.2mm}
$\smash{\Tr(\varrho^{\calT^{(\bbZ_l)}}_{\overline\eta}(\sigma)) = \Tr(\varrho^\calT_{\overline\eta}(\sigma))}$
and may conclude
$\smash{\Tr(\varrho^{\calT^{(\bbZ_l)}}_{\overline\eta}(\sigma)) \not\in \{-22,\ldots,22\} \subset \bbZ_l}$.
Thus,
$\smash{\Tr(\varrho^{\calT^{(\bbZ_l)}}_{\overline\eta}(\sigma))}$
has a positive
\mbox{$l$-adic}
distance from the finite set
$\{-22,\ldots,22\}$,
which means that there exists some
$e\in\bbN$
of the kind that
\begin{equation}
\label{tr_0}
\Tr(\varrho_e(\sigma)) \not\in \{\overline{-22}, \ldots, \overline{22}\} \subset \bbZ/l^e\bbZ \,,
\end{equation}
for 
$\smash{\varrho_e\colon \pi_1^\et(\underline{B}, \overline\eta) \to \GL(\calT^{(\bbZ_l)}_{\overline\eta}/l^e \calT^{(\bbZ_l)}_{\overline\eta})}$
the quotient representation.

We~now apply Proposition~\ref{Cebotarev_Serre}.b) to the product representation
$$\varrho_e \times \chi_K \colon \pi_1^\et(\underline{B}, .) \to \GL(\calT^{(\bbZ_l)}_{\overline\eta}/l^e \calT^{(\bbZ_l)}_{\overline\eta}) \times \Gal(K/\bbQ) \,.$$
Note~that
$\underline{B}$
is certainly flat
over~$\bbZ$,
since it is irreducible and has a
\mbox{$\bbQ$-rational}
point. Proposition~\ref{Cebotarev_Serre}.b) yields a prime
number~$p$
and a closed point
$\tau\colon \Spec\bbF_{\!p} \to \underline{B}$
having the property~that
\begin{equation}
\label{double_congr}
\Frob_\tau \equiv \sigma \pmod {\ker(\varrho_e \!\times\! \chi_K)} \,.
\end{equation}
Here,
$\Frob_\tau$
is to be understood as a suitable representative up to~conjugation.

Formula~(\ref{double_congr}) shows, in particular,
$\Frob_\tau \equiv \sigma \pmod {\ker \chi_K}$,
which implies that
$\Frob_p \in c$.
Moreover,~the congruence modulo
$\ker \varrho_e$
yields
$$\Tr(\varrho_e(\sigma)) = \Tr(\varrho_e(\Frob_\tau)) \,.$$
However,
$\smash{\Tr(\varrho_e(\Frob_\tau)) = \big(\! \Tr(\varrho^{\calT^{(\bbZ_l)}}_{\overline\eta}(\Frob_\tau)) \bmod {l^e} \big) \in \{-\overline{22},\ldots,\overline{22}\} \subset \bbZ/l^e\bbZ}$,
according to the first~step. This~is in contradiction with~(\ref{tr_0}), and the claim is therefore established.\medskip\pagebreak[3]

\noindent
{\em Fourth step.}
Conclusion.

\noindent
In order to complete the argument, we first observe that
$\smash{\GO^0(\calT_{\overline\eta})}$,
equipped with the Zariski topology, is an irreducible topological~space. Indeed, $\smash{\GO^0(\calT_{\overline\eta})}$
is Zariski dense in
$\smash{\GO^0(\calT_{\overline\eta} \otimes_{\bbQ_l}\! \overline\bbQ_l)}$
by~\cite[Corollary in Section~3]{Ro} or \cite[Corollary~2 of Theorem~1]{Che}
and
$\GO^0(\calT_{\overline\eta} \otimes_{\bbQ_l}\! \overline\bbQ_l) \cong \GO^0_{\rk\calT}(\overline\bbQ_l)$
is~irreducible. Consequently,
$\smash{\GO(\calT_{\overline\eta})}$
has at most two components,
$\smash{\GO(\calT_{\overline\eta}) = \dot{\bigcup}{}_{i=1}^N x_i \!\cdot\! \GO^0(\calT_{\overline\eta})}$,
for
$N=1$
or~$2$.

On~the other hand, the result of the previous step implies that
$\smash{\varrho^\calT_{\overline\eta}(U^{(c)}) \subset \GO(\calT_{\overline\eta})}$
cannot be not Zariski dense in any of the components,
$$\overline{\varrho^\calT_{\overline\eta}(U^{(c)})} \cap x_i \!\cdot\! \GO^0(\calT_{\overline\eta}) \subsetneqq x_i \!\cdot\! \GO^0(\calT_{\overline\eta}) \,,$$
for~$i=1,\ldots,N$.
Indeed,
$\diag(\frac\lambda{\rk\calT}, \ldots, \frac\lambda{\rk\calT}) \in \GO_{\rk\calT}(\bbQ_l)$
has trace
$\lambda$,
for an arbitrary
$\lambda\in\bbQ_l$.
Similarly,~when
$\rk\calT >3$
is even,
$$\textstyle \diag(-\frac\lambda{\rk\calT-2}, \frac\lambda{\rk\calT-2}, \ldots, \frac\lambda{\rk\calT-2}) \in \GO_{\rk\calT}(\bbQ_l) \!\setminus\! \GO^0_{\rk\calT}(\bbQ_l)$$
has
trace~$\lambda$.

Consequently,
$\smash{\overline{\varrho^\calT_{\overline\eta}(U)} \cap x_i \!\cdot\! \GO^0(\calT_{\overline\eta}) \subsetneqq x_i \!\cdot\! \GO^0(\calT_{\overline\eta})}$
too, since
$\smash{\overline{\varrho^\calT_{\overline\eta}(U^{(c)})}}$
is a coset of
$\smash{\overline{\varrho^\calT_{\overline\eta}(U)}}$.
Another~application of the same argument shows~that
$$\smash{x \!\cdot\! \overline{\varrho^\calT_{\overline\eta}(U)} \cap \GO^0(\calT_{\overline\eta}) \subsetneqq \GO^0(\calT_{\overline\eta}) \,,}$$
for
every~$\smash{x \in \GO(\calT_{\overline\eta})}$.

Writing
$\smash{\pi_1^\et(\underline{B}, \overline\eta) = \bigcup_{i=1}^{\#\Gal(K/\bbQ)} \!\sigma_i U}$,
one now~finds
$$\MG_{\calT, \underline{B}, l} = \!\!\!\bigcup_{i=1}^{\#\Gal(K/\bbQ)} \!\!\!\!\!\!\!\overline{\varrho^\calT_{\overline\eta}(\sigma_i U)} = \!\!\!\bigcup_{i=1}^{\#\Gal(K/\bbQ)} \!\!\!\!\!\!\!\varrho^\calT_{\overline\eta}(\sigma_i) \!\cdot\! \overline{\varrho^\calT_{\overline\eta}(U)} \,.$$
Therefore,
$$\MG_{\calT, \underline{B}, l} \cap \GO^0(\calT_{\overline\eta}) = \!\!\!\bigcup_{i=1}^{\#\Gal(K/\bbQ)} \!\!\!\!\!\!\!\!\big( \varrho^\calT_{\overline\eta}(\sigma_i) \!\cdot\! \overline{\varrho^\calT_{\overline\eta}(U)} \cap \GO^0(\calT_{\overline\eta}) \big)$$
is the union of finitely many sets, each of which is Zariski closed and properly contained in
$\smash{\GO^0(\calT_{\overline\eta})}$.
Since~$\smash{\GO^0(\calT_{\overline\eta})}$
is an irreducible topological space, this implies
$\MG_{\calT, \underline{B}, l} \cap \GO^0(\calT_{\overline\eta}) \subsetneqq \GO^0(\calT_{\overline{\eta}})$,
which completes the~proof.
\eop
\end{proo1}

\section{The endomorphism field under specialisation}
\label{sec_spec}

It turns out that formula~(\ref{monodromygrp_centr}) from Example~\ref{alg_mon_grp}.v) may be~reversed.

\begin{theo}[\smash{\'Etale cohomological description of the endomorphism field}]
\label{Zarhin_rev}
Let\/~$k$
be a field that is finitely generated
over\/~$\bbQ$,
$X$
a\/
$K3$~surface
defined
over\/~$k$,
and\/
$\calT \subset H^2_\et(X_{\overline{k}}, \bbQ_l(1))$
the transcendental~part. Then
$$C_{\End(\calT)}(\MG_{\calT, k, l}^0) = E \!\otimes_\bbQ\! \bbQ_l \,,$$
for\/
$E$
the endomorphism field of\/
$X(\bbC) = (X \!\times_{\Spec k} \Spec\bbC)(\bbC)$.
Here,
$k$
is arbitrarily embedded
into\/~$\bbC$.
\end{theo}

\begin{coro}[Independence of
$E$
of the embedding into
$\bbC$--elementary case]
\label{E_indep_el}
In~the situation of Theorem~\ref{Zarhin_rev}, let\/
$i_1,i_2\colon k \hookrightarrow \bbC$
be two embeddings and denote by\/
$E_1, E_2$
the corresponding endomorphism~fields.

\begin{abc}
\item
Then\/
$E_1$
and\/~$E_2$
are arithmetically~equivalent.
\item
In~particular,
$[E_1\!:\!\bbQ] = [E_2\!:\!\bbQ]$.
Furthermore,
$E_1$
and\/~$E_2$
have the same normal closure.
If\/~$[E_1\!:\!\bbQ] < 7$
then\/~$E_1 \cong E_2$.
\end{abc}\smallskip

\noindent
{\bf Proof.}
{\em
a)
The isomorphy
\mbox{$E_1 \!\otimes_\bbQ\! \bbQ_l \cong E_2 \!\otimes_\bbQ\! \bbQ_l$},
for every prime
number~$l$,
implies that
$E_1$
and~$E_2$
have the same Dedekind zeta~function. This~is what is called arith\-metic equivalence~\cite{Pe}.\smallskip

\noindent
b)
The~two consequences of arithmetic equivalence are shown in~\cite[Theorem~1]{Pe} and the final statement is \cite[Theorem~1]{BdS}.
}
\eop
\end{coro}

\begin{coro}
\label{Zarhin_E}
In~the situation of Theorem~\ref{Zarhin_rev}, choose, in addition, an embedding
$k \hookrightarrow \bbC$.
Then
$$C_{\End(\calT)}(\MG_{\calT, k, l}^0) \cap \End_\bbQ(T) = E \,. \smallskip$$
\noindent
{\bf Proof.}
{\em
``$\supseteq$'':
Clearly, on one hand, one has
$\smash{E = \End_\Hg(T) \subseteq \End_\bbQ(T)}$
and, on the other,
$\smash{E \subset E \!\otimes_\bbQ\! \bbQ_l = C_{\End(\calT)}(\MG_{\calT, k, l}^0)}$.\smallskip

\noindent
``$\subseteq$'':
Let us put
$\smash{E' := C_{\End(\calT)}(\MG_{\calT, k, l}^0) \cap \End_\bbQ(T)}$.
Then
$\smash{E' \!\otimes_\bbQ\! \bbQ_l}$
is contained in
$\smash{\End_\bbQ(T) \!\otimes_\bbQ\! \bbQ_l = \End(\calT)}$.
Note~here that
$\bbQ_l(1)$
is free of
rank~$1$
over~$\bbQ_l$,
and that we use the identification
$\smash{\calT \cong T \!\otimes_\bbQ\! \bbQ_l(1)}$.
Moreover,~$\smash{E' \!\otimes_\bbQ\! \bbQ_l}$
commutes with
$\smash{\MG_{\calT, k, l}^0}$,
simply because
$E'$
does~so. Therefore,
$$E' \!\otimes_\bbQ\! \bbQ_l \subseteq C_{\End(\calT)}(\MG_{\calT, k, l}^0) = E \!\otimes_\bbQ\! \bbQ_l \,,$$
in view of the~Theorem.
As~$\smash{E,E' \subseteq \End_\bbQ(T)}$
and
$\bbQ_l$
is faithfully flat
over~$\bbQ$,
this yields that
$E' \subseteq E$.%
}%
\eop
\end{coro}

Recall that an extension field
$k$
of~$\bbQ$
is called {\em primary,} if it does not contain any proper algebraic extension
of~$\bbQ$
\cite[\S4.3]{EGAIV}.

\begin{coro}[Independence of
$E$
of the embedding into
$\bbC$--primary
case]
\label{E_indep}
Let\/~$k$
be a field that is finitely generated and primary
over\/~$\bbQ$
and\/
$X$
a\/
$K3$~surface
over\/~$k$.
Then~the endomorphism field of\/
$X(\bbC) = (X \times_{\Spec k} \Spec\bbC)(\bbC)$
is independent of the embedding\/
$k \hookrightarrow \bbC$~chosen.\medskip

\noindent
{\bf Proof.}
{\em
Take~an integral scheme
$B$
of finite type
over~$\bbQ$
with function
field~$k$.
Then,~as
$k$~is
primary, \cite[Proposition~4.5.9]{EGAIV} shows that
$B$
is geometrically~irre\-ducible. Moreover,~according to H.\ Hironaka~\cite{Hi}, one may assume that
$B$~is
non\-sin\-gu\-lar.

The~usual spreading out argument provides a morphism
$q\colon \calX \to B$
of
\mbox{$\bbQ$-schemes}
of finite type with generic
fibre~$X$.
Restricting~$B$
to an open subscheme, if necessary, one may assume that
$q$~is
proper and~smooth and that every fibre is a
$K3$~surface.

Next,~consider two embeddings
$i_1, i_2\colon k = k(\eta) \hookrightarrow \bbC$.
These~yield two complex points
$\smash{\eta_1^c, \eta_2^c\colon \Spec\bbC \to B}$
on~$B$,
and hence on
$\smash{B_\bbC := B \times_{\Spec\bbQ} \Spec\bbC}$,
as well as on the complex manifold
$\smash{B(\bbC) = B_\bbC(\bbC)}$.
As~$B(\bbC)$
is connected, we may choose a path
$w \in \pi_1(B(\bbC), \eta_1^c, \eta_2^c)$.

The higher direct image sheaf
$R^2 q(\bbC)_* \bbQ$
is locally free
on~$B(\bbC)$
of rank~$22$.
Thus,~the
path~$w$
induces an
isomorphism~$i_w$
between the stalks
$\smash{(R^2 q(\bbC)_* \bbQ)_{\eta_1^c}}$
and
$\smash{(R^2 q(\bbC)_* \bbQ)_{\eta_2^c}}$,
which are, according to Grauert's Theorem \cite[Satz~5]{Gr}, canonically isomorphic to
$\smash{H^2(\calX_{\eta_1^c}(\bbC), \bbQ)}$
and
$\smash{H^2(\calX_{\eta_2^c}(\bbC), \bbQ)}$,~respectively.

On~the other hand, as usual,
$w$
induces an \'etale path on
$\smash{B_\bbC}$,
and, via the natural projection, an element
$\smash{\widetilde{w} \in \pi_1^\et(B, \overline\eta)}$.
By~\cite[Expos\'e~V, Proposition~8.2]{SGA1}, this fundamental group is a quotient of
$\smash{\Gal(\overline{k}/k)}$.
Thus,~$\smash{\widetilde{w}}$
operates as an automorphism
of~$\smash{\overline{k}}$
and hence preserves the algebraic classes
in~$\smash{H^2_\et(\calX_{\overline\eta}, \bbQ_l(1)) = H^2_\et(X_{\overline{k}}, \bbQ_l(1))}$.

Consequently, one also has that
$\smash{i_{\widetilde{w}} (\calT) = \calT}$.
I.e., the isomorphism
$$\smash{i_w \colon H^2(X_{\eta_1^c}(\bbC), \bbQ) \to H^2(X_{\eta_2^c}(\bbC), \bbQ)}$$
maps
$\smash{T_{\eta_1^c} \!\otimes_\bbQ\! \bbQ_l(1)}$
onto
$\smash{T_{\eta_2^c} \!\otimes_\bbQ\! \bbQ_l(1)}$.
As~$\bbQ_l(1)$
is faithfully flat over
$\bbQ$,
this implies
$\smash{i_w (T_{\eta_1^c}) = T_{\eta_2^c}}$.

Conjugation~by
$i_w$
hence provides an isomorphism
$$c_w \colon \End_\bbQ(T_{\eta_1^c}) \to \End_\bbQ(T_{\eta_2^c}) \,, \quad M \mapsto i_w \!\circ\! M \!\circ\! i_w^{-1} \,,$$
which clearly induces the isomorphism
$\smash{c_{\widetilde{w}} \colon \End(\calT) \to \End(\calT), \, M \mapsto i_{\widetilde{w}} \!\circ\! M \!\circ\! i_{\widetilde{w}}^{-1}}$.
Moreover,~as noticed in Remark~\ref{indep_basepoint},
$c_{\widetilde{w}}$
maps the algebraic monodromy group
$\MG_{\calT, k, l}$,
and therefore also the identity component
$\MG_{\calT, k, l}^0$,
onto~itself. Consequently,
$$c_w \big( C_{\End(\calT)}(\MG_{\calT, k, l}^0) \cap \End_\bbQ(T_{\eta_1^c}) \big) \subseteq C_{\End(\calT)}(\MG_{\calT, k, l}^0) \cap \End_\bbQ(T_{\eta_2^c}) \,.$$
Corollary~\ref{Zarhin_E} shows that this means nothing but
$c_w(E_1) \subseteq E_2$.
As~$c_w$
is injective and one has
$[E_1 \!:\! \bbQ] = [E_2 \!:\! \bbQ]$,
due to Corollary~\ref{E_indep_el}.b), the assertion~follows.
}%
\eop
\end{coro}

In order to draw more conclusions from Theorem~\ref{Zarhin_rev}, we need a~lemma.

\begin{lem}
\label{tr_jump}
Let\/~$k$
be a field that is finitely generated
over\/~$\bbQ$,
$X$
a\/
$K3$~surface
over\/~$k$,
and\/
$\calT' \subseteq H^2_\et(X_{\overline{k}}, \bbQ_l(1))$
a subvector space on which the cup product pairing is non-degenerate and which contains the transcendental
part\/~$\calT$.
Let\/
$\calR \subset \calT'$
be the orthogonal complement
of\/
$\calT$
in\/~$\calT'$.

\begin{abc}
\item
Then\/
$\smash{\MG_{\calT', k, l}^0 = \MG_{\calT, k, l}^0}$,
the operation of\/
$\smash{\MG_{\calT', k, l}^0}$
on\/~$\calR$
being~trivial.
\item
The centraliser\/
$\smash{C_{\End(\calT')}(\MG_{\calT', k, l}^0)}$
maps
$\calR$
to~itself.
\end{abc}\smallskip

\noindent
{\bf Proof.}
{\em
a)
Being~perpendicular to
$\calT$,
the direct summand
$\calR$
consists of algebraic~classes. These~are defined over a finite extension field
of~$k$,
and therefore pointwise fixed under an open subgroup of finite index
in~$\smash{\MG_{\calT', k, l}}$.
In~other words,
$\MG_{\calT', k, l}^0$
operates as the identity map
on~$\calR$.
This~yields
$\smash{\MG_{\calT', k, l}^0 = \MG_{\calT, k, l}^0}$,
since the operation of Galois preserves orthogonality.\smallskip

\noindent
b)
On~the other hand, no non-zero element of
$\calT$
is fixed under
$\smash{\MG_{\calT', k, l}^0}$.
Indeed,~such an element
$x \in \calT$
would be fixed under
$\smash{\Gal(\overline{k}/k')}$,
for
$k'$
a certain finite extension field
of~$k$,
and hence algebraic, according to the Tate conjecture (cf.\ Fact~\ref{Tate}.a)), a contradiction.
As~elements commuting with
$\smash{\MG_{\calT', k, l}^0}$
cannot interchange fixed points with non-fixed points, we see that
$\smash{C_{\End(\calT')}(\MG_{\calT', k, l}^0)}$
maps
$\calR$
to~itself.
}
\eop
\end{lem}

\begin{coro}[The endomorphism field under specialisation]
\label{endfield_spec}
Let\/~$q\colon X \to B$
be a proper and smooth morphism of geometrically connected schemes of finite type over\/~$\bbQ$,
every fibre of which is a\/
$K3$~surface.

\begin{abc}
\item
Then~the endomorphism field\/
$E$
of the generic fibre
$X_\eta$
is independent of the embedding\/
$k(\eta) \hookrightarrow \bbC$.
\item
Let\/
$s \in B$
be a~point. Choose~arbitrary embeddings\/
$k(\eta) \hookrightarrow \bbC$
and\/
$k(s) \hookrightarrow \bbC$
and denote~by\/
$\eta^c, s^c \in B(\bbC)$
the complex points corresponding to\/
$\eta$
and\/~$s$,
respectively. Let,~moreover\/
$w \in \pi_1(B(\bbC), \eta^c, s^c)$
be a~path.
\begin{iii}
\item
Then\/~$w$
induces an isomorphism\/
$\smash{i_w\colon H^2(X_{\eta^c}(\bbC), \bbQ) \to H^2(X_{s^c}(\bbC), \bbQ)}$,
which maps the transcendental part\/
$\smash{T^{(X_{\eta^c})} \subset H^2(X_{\eta^c}(\bbC), \bbQ)}$
to some\/
$\smash{T_{s^c} \supseteq T^{(X_{s^c})}}$.
\item
Thus, by transport of structure,
$E$
operates
on\/~$\smash{T_{s^c}}$,
too. Under~this operation,
$\smash{T^{(X_{s^c})}}$
is mapped to~itself.
\item
For~the endomorphism field\/
$\smash{E^{(X_{s^c})}}$
of\/~$X_s(\bbC)$,
one has\/
$\smash{E^{(X_{s^c})} \supseteq E}$.
Thereby,~the operation
of\/~$\smash{E \subseteq E^{(X_{s^c})}}$
coincides with that obtained by transport of~structure.
\end{iii}
\end{abc}\smallskip

\noindent
{\bf Proof.}
{\em
a)
As~$B$
is geometrically connected,
$k(\eta)$
is clearly primary
over~$\bbQ$.
Thus, Corollary~\ref{E_indep} implies the~assertion.\smallskip

\noindent
b.i)
The higher direct image sheaf
$R^2 q(\bbC)_* \bbQ$
is locally free
on~$B(\bbC)$
of
rank~$22$.
Thus,~$w$
induces an isomorphism
$i_w\colon (R^2 q(\bbC)_* \bbQ)_{\eta^c} \to (R^2 q(\bbC)_* \bbQ)_{s^c}$.
Note~that these stalks are canonically isomorphic to
$H^2(X_{\eta^c}(\bbC), \bbQ)$
and
$H^2(X_{s^c}(\bbC), \bbQ)$,
respectively, due to \cite[Satz~5]{Gr}.

Moreover, algebraic classes remain algebraic under specialisation, i.e.\
\begin{equation}
\label{Pic_spec}
i_w (P^{(X_{\eta^c})}) \subseteq P^{(X_{s^c})} \,,
\end{equation}
in the situation of a smooth~family. Indeed, they are representable by Weil divisors and one may just take the Zariski closure of a representing Weil~divisor. The~claim follows immediately from (\ref{Pic_spec}) by taking orthogonal complements on both~sides.\smallskip

\noindent
ii)
The~assertion descends under base change by a finite morphism
$\smash{p\colon B' \to B}$.
Thus,~we may assume without restriction that
$B$
is a normal scheme~\cite[Chapitre~V, \S1, Corollaire~1 du Proposition~18]{Bo}.

Let~us switch to \'etale cohomology. The~algebraic part
$\smash{\calP_{\overline\eta} \subset H^2_\et(X_{\overline\eta}, \bbQ_l(1))}$
extends to a locally constant
sheaf~$\calP$
on the whole
of~$B$,
by virtue of of Lemma~\ref{sheaf_ext}. Consequently,~the transcendental part
$\smash{\calT_{\overline\eta} \subset H^2_\et(X_{\overline\eta}, \bbQ_l(1))}$
extends to a twisted-constant
sheaf~$\calT$.
Moreover, the comparison theorem between \'etale and complex cohomology \cite[Expos\'e~11, Th\'eor\`eme~4.4.iii)]{SGA4} shows that the commutative diagram of the given data
$$
\xymatrixcolsep{2.7mm}
\xymatrixrowsep{3.0mm}
\xymatrix{
H^2(X_{\eta^c}(\bbC), \bbQ) \ar@{->}[rr]^{i_w} && H^2(X_{s^c}(\bbC), \bbQ) \\
T^{(X_{\eta^c})} \ar@{->}[rr]^{i_w} \ar@{^{(}->}[u] && T_{s^c} \ar@{^{(}->}[u] \\
&& T^{(X_{s^c})} \ar@{^{(}->}[u] \phantom{\,.}
}
$$
goes over under tensoring
with~$\bbQ_l(1)$
into
$$
\xymatrixcolsep{2.7mm}
\xymatrixrowsep{3.0mm}
\xymatrix{
H^2_\et(X_{\overline\eta}, \bbQ_l(1)) \ar@{->}[rr] && H^2_\et(X_{\overline{s}}, \bbQ_l(1)) \\
\calT_{\overline\eta} \ar@{->}[rr] \ar@{^{(}->}[u] && \calT_{\overline{s}} \ar@{^{(}->}[u] \\
&& \calT^{(X_{\overline{s}})} \ar@{^{(}->}[u] \,.
}
$$

On~the other hand, we have
$\smash{\MG_{\calT_{\overline{s}}, k(s), l}^0 \subseteq \MG_{\calT, B, l}^0 = \MG_{\calT_{\overline\eta}, k(\eta), l}^0}$,
due to~\ref{bc_spec}.i) and~ii). Here, in view of Remark \ref{indep_basepoint}, the underlying identification may be supposed to be induced by
$w$.
The~inclusion~yields
\begin{equation}
\label{cent_spec}
C_{\End(\calT_{\overline{s}})}(\MG_{\calT_{\overline{s}}, k(s), l}^0) \supseteq C_{\End(\calT_{\overline\eta})}(\MG_{\calT_{\overline\eta}, k(\eta), l}^0) = E \!\otimes_\bbQ\! \bbQ_l \,.
\end{equation}

Now~write
$\smash{\calT_{\overline{s}} = \calT^{(X_{\overline{s}})} \oplus \calR}$,
with a direct summand
$\calR$
that is perpendicular
to~$\smash{\calT^{(X_{\overline{s}})}}$.
Then,~according to Lemma~\ref{tr_jump},
$\smash{C_{\End(\calT_{\overline{s}})}(\MG_{\calT_{\overline{s}}, k(s), l}^0)}$
maps
$\calR$
to~itself.
Since~$\smash{E \!\otimes_\bbQ\! \bbQ_l}$
acts via self-adjoint endomorphisms~\cite[Theorem~1.5.1]{Za}, this shows that
$\smash{E \!\otimes_\bbQ\! \bbQ_l}$
maps
$\smash{\calT^{(X_{\overline{s}})}}$
to itself,~either.

Translating~this back to complex cohomology, one finds that
$\smash{E \subseteq \End(T_{s^c})}$
has the property that
$\smash{E \subset E \!\otimes_\bbQ\! \bbQ_l \subseteq \End(T_{s^c} \!\otimes_\bbQ\! \bbQ_l(1))}$
maps
$\smash{T^{(X_{s^c})} \!\otimes_\bbQ\! \bbQ_l(1)}$
into~itself.
As~$\bbQ_l(1)$
is faithfully flat over
$\bbQ$,
this is enough to enforce the claim.\smallskip

\noindent
iii)
By~(\ref{cent_spec}), the operation of
$\smash{E \subset E \!\otimes_\bbQ\! \bbQ_l}$
on~$\smash{\calT_{\overline{s}}}$
commutes with
$\smash{\MG_{\calT_{\overline{s}}, k(s), l}^0}$.
Moreover,~according to Lemma~\ref{tr_jump}.a),
$\smash{\MG_{\calT_{\overline{s}}, k(s), l}^0}$
maps
$\smash{\calT^{(X_{\overline{s}})} \subseteq \calT_{\overline{s}}}$
to~itself, while
$E \otimes_\bbQ\! \bbQ_l$
does the same, as shown in~b). Restricting the endomorphisms
to~$\smash{\calT^{(X_{\overline{s}})}}$,
we find that
$\smash{E \subset E \!\otimes_\bbQ\! \bbQ_l \subseteq \End(\calT^{(X_{\overline{s}})})}$
commutes with
$\smash{\MG_{\calT^{(X_{\overline{s}})}, k(s), l}^0}$.
In~other words,
$$\smash{E \subset E \!\otimes_\bbQ\! \bbQ_l \subseteq C_{\End(\calT^{(X_{\overline{s}})})} \big(\! \MG_{\calT^{(X_{\overline{s}})}, k(s), l}^0 \!\big) \,.}$$
Since,~according to~b),
$E$~maps
$\smash{T^{(X_{s^c})}}$
to itself, this~yields
$$E \subseteq C_{\End(\calT^{(X_{\overline{s}})})} \big(\! \MG_{\calT^{(X_{\overline{s}})}, k(s), l}^0 \!\big) \cap \End_\bbQ(T^{(X_{s^c})}) = E^{(X_{s^c})} \,,$$
as~claimed.
}%
\eop
\end{coro}

\begin{coro}[Complex fibres]
\label{complex}
Let\/~$q\colon X \to B$
be a proper and smooth morphism of geometrically connected schemes of finite type
over\/~$\bbQ$,
every fibre of which is a\/
$K3$~surface.
Suppose~that the generic fibre\/
$X_\eta$
has real or complex multiplication by an endomorphism
field\/~$\smash{E \supsetneqq \bbQ}$.
Then,~for every complex point\/
$\theta\in B(\bbC)$,
the fibre\/
$X_\theta(\bbC)$
is acted upon
by\/~$E$.\medskip

\noindent
{\bf Proof.}
{\em
We~suppose that
$\dim B \geq 1$,
as otherwise there is nothing to~prove. Take~an open neighbourhood
$U \ni \theta$
that is connected and simply~connected. Then
$\smash{R^2 (q |_{q^{-1}(U)})_* \bbQ = (\bbQ^{22})_U}$
is a constant~sheaf. 

Moreover,~the subset
$\smash{U^\alg := B(\overline\bbQ) \cap U}$
of algebraic points is dense
in~$U$
with respect to the complex topology.
$U^\alg$~is
the same as the set of points of
type~$s^c$,
for~$s \in B$
a closed~point. For~each
$z \in U^\alg$,
the fibre
$X_z(\bbC)$
is acted upon
by~$E$,
as shown in Corollary~\ref{endfield_spec}.

In~order to make this more precise, let us take an embedding
$k(\eta) \hookrightarrow \bbC$
and denote by
$\eta^c \in B(\bbC)$
the corresponding complex~point. In~addition, we choose one particular
point~$z_0 \in U^\alg$
and a path
$w \in \pi_1(B(\bbC), \eta^c, z_0)$.
For~every other point
$z \in U^\alg$,
we choose a path
$w_z \in \pi_1(U, z_0, z)$.
Let~us denote the element
of~$\pi_1(B(\bbC), z_0, z)$,
induced
by~$w_z$,
again
by~$w_z$.
Then~Corollary~\ref{endfield_spec} applies to
$w_z \!\circ\! w \in \pi_1(B(\bbC), \eta^c, z)$.

It~shows that the fibre
$X_z(\bbC)$
is acted upon
by~$E$
in the following~manner.
One~has
$i_{w_z \circ w} (T^{(X_{\eta^c})}) \supseteq T^{(X_z)}$
and the operation
of~$E$
on
$i_{w_z \circ w} (T^{(X_{\eta^c})})$,
induced by that
on~$T^{(X_{\eta^c})}$,
defines the action
on~$T^{(X_z)}$.
Thus,~the actions
of~$E$
on all
$i_{w_z \circ w} (T^{(X_{\eta^c})})$,
for
$z \in U^\alg$,
are compatible among each other, via transport of structure
under~$w_z$.
But~the latter is the obvious one on the constant sheaf
$\smash{R^2 (q |_{q^{-1}(U)})_* \bbQ}$.

The~operation
of~$E$
splits
$\smash{i_w (T^{(X_{\eta^c})}) \!\otimes_\bbQ\! \bbC}$
into
$r = [E\!:\!\bbQ]$
eigenspaces
$V_1,\ldots,V_r$.
The~same decomposition applies to every
$z \in U^\alg$,
$$i_{w \circ w_z} (T^{(X_{\eta^c})}) \!\otimes_\bbQ\! \bbC = \bigoplus_{i=1}^r ((V_i)_U)_z \,.$$
As~$\smash{E = \End_\Hg(T^{(X_z)})}$,
this means that
$H^{2,0}(X_z(\bbC),\bbC) \subseteq ((V_i)_U)_z$,
for a
certain~$i$.
In~other words, the one-dimensional vector space
$$H^{2,0}(X_z(\bbC),\bbC) \in \Pb(H^2(X_z(\bbC),\bbC)) \cong \Pb(H^2(X_{z_0}(\bbC),\bbC)) \,,$$
represents a point lying on the union of
the~$r$
projective subspaces
$\Pb(V_1), \ldots, \Pb(V_r)$.

On~the other hand, the mapping
$\Pi\colon U \to \Pb(H^2(X_{z_0}(\bbC),\bbC)), z \mapsto\! H^{2,0}(X_z(\bbC),\bbC)$,
is holomorphic \cite[Theorem~IV.4.2]{BHPV}. Our~argument will only need~continuity. Indeed,~we have
$\Pi(U^\alg) \subseteq \bigcup_{i=1}^r \Pb(V_i)$.
As~the right hand side is a closed subset and
$\Pi$~is
continuous, this implies that
$\Pi(U) \subseteq \bigcup_{i=1}^r \Pb(V_i)$.
But~this means that
$X_z(\bbC)$
is acted upon
by~$E$,
for any
$z \in U$,
and in particular
for~$z = \theta$.
}
\eop
\end{coro}

In~order to prove Theorem~\ref{Zarhin_rev}, we need a few auxiliary~results.

\begin{subl}
\label{span_On}
One has\/
$\smash{\spann_{\overline\bbQ_l} \SO_n(\overline\bbQ_l) = \M_{n\times n}(\overline\bbQ_l)}$,
for every natural number\/
$n \neq 2$.\medskip

\noindent
{\bf Proof.}
{\em
Put
$\smash{S := \spann_{\overline\bbQ_l} \SO_n(\overline\bbQ_l)}$.
Then~clearly
$S \subseteq \M_{n\times n}(\overline\bbQ_l)$.
Moreover,
$S$
is not just a
\mbox{$\overline\bbQ_l$-vector}
space, but a representation of
$\SO_n(\overline\bbQ_l) \!\times\! \SO_n(\overline\bbQ_l)$,
via
$$(\SO_n(\overline\bbQ_l) \!\times\! \SO_n(\overline\bbQ_l)) \times S \longrightarrow S\,, \quad ((M_1,M_2), s) \mapsto M_1sM_2^{-1} \,.$$
Thus it suffices to show that
$\M_{n\times n}(\overline\bbQ_l)$
is irreducible as an
$\SO_n(\overline\bbQ_l) \!\times\! \SO_n(\overline\bbQ_l)$-rep\-re\-sen\-ta\-tion.
But~$\smash{\M_{n\times n}(\overline\bbQ_l) \cong \overline\bbQ_l^n \!\otimes\! (\overline\bbQ_l^n)^*}$
and
$\smash{\overline\bbQ_l^n}$
is irreducible as a
$\smash{\SO_n(\overline\bbQ_l)}$-module,
due to \cite[Chapter~VI, (5.4.v)]{BtD}, for
$n\ge3$,
and trivially,
for~$n=1$.
}
\eop
\end{subl}

\begin{lem}
\label{span_algabg}
Let\/~$\calT$
be a finite-dimensional\/
\mbox{$\bbQ_l$-vector}
space equipped with a non-degenerate symmetric form that is acted upon by a totally real or CM field\/
$E$
via self-adjoint linear~maps. In~the case that\/
$E$
is totally real, suppose that\/
$\calT$
is free as an\/
\mbox{$E \!\otimes_\bbQ\! \bbQ_l$-module}
of
rank\/~$\neq\!2$.
Then
$$\smash{\spann_{\overline\bbQ_l} (C_{\O(\calT \otimes_{\bbQ_l} \overline\bbQ_l)}(E))^0 = C_{\End(\calT \otimes_{\bbQ_l} \overline\bbQ_l)}(E) \,.}$$.

\noindent
{\bf Proof.}
{\em
{\em First case:}
$E$
is totally~real.

\noindent
Then
$\smash{\calT \!\otimes_{\bbQ_l}\! \overline\bbQ_l}$
is split under the operation
of~$E$
into
$r = [E:\bbQ]$
simultaneous eigenspaces
$\smash{V_1,\ldots,V_r \subset \calT \!\otimes_{\bbQ_l}\! \overline\bbQ_l}$,
which are mutually perpendicular, due to the self-adjointness assumption. Hence,
$$C_{\End(\calT \otimes_{\bbQ_l} \overline\bbQ_l)}(E) = \{f\colon \!\calT \!\otimes_{\bbQ_l}\! \overline\bbQ_l \!\to\! \calT \!\otimes_{\bbQ_l}\! \overline\bbQ_l \mid f \;\overline\bbQ_l\text{-linear}, f(V_i) \subseteq V_i \,\text{for}\, i=1,\ldots,r\}$$
and
$\smash{(C_{\O(\calT \otimes_{\bbQ_l} \overline\bbQ_l)}(E))^0}$
is given analogously, with the additional assumption that
all restrictions
$f |_{V_i}\colon V_i \to V_i$
be orthogonal maps of
determinant~$1$.
Since\/
$\dim V_i = \rk_{E \otimes_\bbQ \overline\bbQ_l} (\calT \!\otimes_{\bbQ_l}\! \overline\bbQ_l) = \rk_{E \otimes_\bbQ \bbQ_l} \calT \neq 2$,
Sublemma~\ref{span_On} implies the~claim.\smallskip\pagebreak[3]

\noindent
{\em Second case:}
$E$
is a CM~field.

\noindent
Here,
$\smash{\calT \!\otimes_{\bbQ_l}\! \overline\bbQ_l}$
is split under the operation
of~$E$
into
$2s = [E:\bbQ]$
simultaneous eigenspaces
$\smash{V_1,\overline{V}_1,\ldots,V_s,\overline{V}_s \subset \calT \!\otimes_{\bbQ_l}\! \overline\bbQ_l}$.
These~are isotropic and mutually perpendicular, with the only exceptions that
$V_i \not\perp \overline{V}_i$,
for~$i = 1,\ldots,s$.
Indeed,~if some primitive element
$e \in E$
acts with eigenvalues
$\lambda_i$
and~$\lambda_j$
on
$V_i$
and~$V_j$,
respectively, then one~finds
$$\lambda_i \langle v_i, v_j\rangle = \langle e v_i, v_j\rangle = \langle v_i, e v_j\rangle = \overline\lambda_j \langle v_i, v_j\rangle \,,$$
which yields
$\langle v_i, v_j\rangle = 0$,
as soon as
$\lambda_i \neq \overline\lambda_j$.

Hence,
\begin{align*}
(C_{\O(\calT \otimes_{\bbQ_l} \overline\bbQ_l)}(E))^0\hspace{3.5cm} & \\[-2mm]
= \{f\colon \!\calT \!\otimes_{\bbQ_l}\! \overline\bbQ_l \!\to\! \calT \!\otimes_{\bbQ_l}\! \overline\bbQ_l \mid{} & f \;\overline\bbQ_l\text{-linear}, f(V_i) \subseteq V_i, f(\overline{V}_i) \subseteq \overline{V}_i \\[-2mm]
& \text{and } f |_{V_i \oplus \overline{V}_i} \in \SO(V_i \!\oplus\! \overline{V}_i) \text{ for } i=1,\ldots,r\}
\end{align*}
and
$\smash{C_{\End(\calT \otimes_{\bbQ_l} \overline\bbQ_l)}(E)}$
is given in the same way, dropping the second~condition.

But~all matrices of~type
$\smash{(\atop{A \;\;\;0\;\;\;\;\;\;\;}{0 \,(A^t)^{-1}}\!)}$,
for
$A \in \GL_{\dim V_i}(\overline\bbQ_l)$,
represent orthogonal maps of
determinant~$1$,
as a direct calculation~shows. The~assertion immediately follows from~this.
}
\eop
\end{lem}

\begin{coro}
\label{span_algabgGO}
In the situation of Lemma~\ref{span_algabg}, one~has
$$\spann_{\overline\bbQ_l} (C_{\GO(\calT \otimes_{\bbQ_l} \overline\bbQ_l)}(E))^0 = C_{\End(\calT \otimes_{\bbQ_l} \overline\bbQ_l)}(E) \,. \eop$$
\end{coro}\vspace{-.7cm}

\begin{prop}
\label{span_Ql}
Let\/~$\calT$
be a finite-dimensional\/
\mbox{$\bbQ_l$-vector}
space equipped with a non-degenerate symmetric form that is acted upon by a totally real or CM field\/
$E$
via self-adjoint linear~maps. In~the case that\/
$E$
is totally real, suppose that\/
$\calT$
is free as an\/
\mbox{$E \!\otimes_\bbQ\! \bbQ_l$-module}
of
rank\/~$\neq\!2$.
Then\/
$$\spann_{\bbQ_l} (C_{\GO(\calT)}(E))^0 = C_{\End(\calT)}(E) \,.$$

\noindent
{\bf Proof.}
{\em
The inclusion
``$\subseteq$''
is~obvious. In~order to prove
``$\supseteq$'',
assume the~contrary. Then
$$(C_{\GO(\calT)}(E))^0 \subseteq \spann_{\bbQ_l}\! (C_{\GO(\calT)}(E))^0 \subseteq C_{\End(\calT)}\!(E) \cap V(\lambda) \subseteq C_{\End(\calT \!\otimes_{\bbQ_l}\! \overline\bbQ_l)}\!(E) \cap V(\lambda) \,,$$
for a certain
\mbox{$\bbQ_l$-linear}
form~$\lambda$
on the
\mbox{$\bbQ_l$-vector}
space
$\smash{C_{\End(\calT)}(E)}$.
For~the Zariski closure, this~shows
$$\overline{(C_{\GO(\calT)}(E))^0} \subseteq C_{\End(\calT \otimes_{\bbQ_l} \overline\bbQ_l)}(E) \cap V(\lambda) \subsetneqq C_{\End(\calT \otimes_{\bbQ_l} \overline\bbQ_l)}(E) \,.$$

On~the other hand,
$\smash{(C_{\GO(\calT \otimes_{\bbQ_l} \overline\bbQ_l)}(E))^0}$
is the set of all
\mbox{$\overline\bbQ_l$-rational}
points on a connected linear algebraic group that is defined
over~$\bbQ_l$,
while
$\smash{(C_{\GO(\calT)}(E))^0}$
is the set of all
\mbox{$\bbQ_l$-rational}
points. 
Therefore,~it is known~\cite{Ro,Che} that
$\smash{(C_{\GO(\calT)}(E))^0}$
is Zariski dense
in~$\smash{(C_{\GO(\calT \otimes_{\bbQ_l} \overline\bbQ_l)}(E))^0}$.
Hence,
$$(C_{\GO(\calT \otimes_{\bbQ_l} \overline\bbQ_l)}(E))^0 \subseteq C_{\End(\calT \otimes_{\bbQ_l} \overline\bbQ_l)}(E) \cap V(\lambda) \subsetneqq C_{\End(\calT \otimes_{\bbQ_l} \overline\bbQ_l)}(E) \,,$$
implying
$\smash{\spann_{\overline\bbQ_l} (C_{\GO(\calT \otimes_{\bbQ_l} \overline\bbQ_l)}(E))^0 \subseteq C_{\End(\calT \otimes_{\bbQ_l} \overline\bbQ_l)}(E) \cap V(\lambda) \subsetneqq C_{\End(\calT \otimes_{\bbQ_l} \overline\bbQ_l)}(E)}$\,,
too. Thus,~we arrived at a contradiction with Corollary~\ref{span_algabgGO}.
}
\eop
\end{prop}

\begin{proo2}
By~(\ref{monodromygrp_centr}),
$C_{\End(\calT)}(\MG_{\calT, k, l}^0) = C_{\End(\calT)}((C_{\GO(\calT)}(E))^0)$.
But every endomorphism
of~$\calT$
that commutes with a certain set, commutes with its
\mbox{$\bbQ_l$-span}.
I.e.,
$$C_{\End(\calT)}(\MG_{\calT, k, l}^0) = C_{\End(\calT)}(\spann_{\bbQ_l} (C_{\GO(\calT)}(E))^0) \,.$$
Moreover,
$E$
can only be a totally real or a CM field by \cite[Theorems~1.6.a) and 1.5.1]{Za},
$\calT = T \!\otimes_\bbQ\! \bbQ_l(1)$
is free over
$E \!\otimes_\bbQ\! \bbQ_l$,
and in the totally real case
$\dim_E T = 2$
is not possible \cite[Lemma~3.2]{vG}. Thus, Proposition~\ref{span_Ql} yields~that
$$C_{\End(\calT)}(\MG_{\calT, k, l}^0) = C_{\End(\calT)}(\spann_{\bbQ_l} (C_{\GO(\calT)}(E))^0) = C_{\End(\calT)}(C_{\End(\calT)}(E)) \,.$$
Moreover,~the elements of
$C_{\End(\calT)}(E)$
actually commute with
$E \!\otimes_\bbQ\! \bbQ_l \supset E$,
which is still contained
in~$\End(\calT)$.
I.e.,
$C_{\End(\calT)}(E) = C_{\End(\calT)}(E \!\otimes_\bbQ\! \bbQ_l)$.
Finally,~the classical double centraliser theorem \cite[Corollary~13.18]{Is} applies and shows~that
$$C_{\End(\calT)}(\MG_{\calT, k, l}^0) = C_{\End(\calT)}(C_{\End(\calT)}(E \!\otimes_\bbQ\! \bbQ_l)) = E \!\otimes_\bbQ\! \bbQ_l \,. \eop$$
\end{proo2}\vspace{-2\bigskipamount}

\section{An explicit family}
\label{sec_exfam_Qw2}

The~example below gives an illustration on how to apply Theorem~\ref{suff_crit}. The~family considered has been treated before~\cite{EJ14}, so most of its properties needed in the proof can simply be~cited.

\begin{ex}[An explicit family of
$K3$~surfaces
with RM
by~$\bbQ(\sqrt{2})$]
\label{Qw2_family}
Let
$\smash{q\colon X \!\to\!\! B}$,
for
$\smash{B \!:=\! \Spec \bbQ[t,\!\frac1{t(t^2-2)(t^2+2)(t^2-4t+2)(t^2+4t+2)}] \!\subset\! \Ab^1_\bbQ}$,
be the family of
$K3$~surfaces
that is fibre-by-fibre the minimal desingularisation of the double cover
of~$\Pb^2$,
given~by
\begin{eqnarray*}
w^2 & = & \textstyle
 [(\frac18 t^2 \!-\! \frac12 t \!+\! \frac14)y^2 + (t^2 \!-\! 2t \!+\! 2)yz + (t^2 \!-\! 4t \!+\! 2)z^2] \\[-1mm]
 & & \textstyle\hspace{4mm}
 [(\frac18 t^2 \!+\! \frac12 t \!+\! \frac14)x^2 + (t^2 \!+\! 2t \!+\! 2)xz + (t^2 \!+\! 4t \!+\! 2)z^2] [2x^2 + (t^2 \!+\! 2)xy + t^2y^2] \, .
\end{eqnarray*}
Then

\begin{iii}
\item
the generic fibre
$X_\eta$
of~$q$
is of geometric Picard
rank~$16$.
\item
The endomorphism field of
$X_\eta$
is~$\smash{\bbQ(\sqrt{2})}$.
\item
For~every
$\theta \in B(\bbC)$,
the transcendental part
$\smash{T \subset H^2(X_\theta(\bbC), \bbQ)}$
of the cohomology of the fibre
$X_\theta(\bbC)$
is acted upon
by~$\smash{\bbQ(\sqrt{2})}$.
\item
Let the complex point
$\theta \in B(\bbC)$
be of the kind that the fibre
$X_\theta(\bbC)$
has Picard
rank~$16$.
Then
$X_\theta(\bbC)$
has real multiplication
by~$\smash{\bbQ(\sqrt{2})}$.
\end{iii}\smallskip

\noindent
{\bf Proof.}
For~$t_0 \in B$,
the ramification locus of the double cover
underlying~$X_{t_0}$
is a union of six lines, all of which are defined over
$\smash{\bbQ(\sqrt{2},t_0)}$.
Moreover,~a direct calculation shows that no three of these lines have a point in common, except for
$t_0 = 0$,
$\smash{\pm\sqrt{2}}$,
$\smash{\pm\sqrt{-2}}$,
and
$\smash{\pm2\!\pm\!\sqrt{2}}$,
which are exactly the points excluded
from~$B$.
Thus,~$q$
is indeed a well-defined family of
$K3$~surfaces.\smallskip

\noindent
i)
Certainly,
$\Pic X_{\overline\eta}$
contains the
$\smash{(\atop{6}2) = 15}$
classes of the exceptional curves, obtained by blowing up the intersection points of the ramification locus, and the pull-back of the class of the general line
on~$\Pb^2$.
Hence,~$\rk\Pic X_{\overline\eta} \geq 16$.
On~the other hand, for
$t_0 \equiv 1 \pmod {17\!\cdot\!23}$,
one has
$\smash{\rk\Pic X_{{\overline{t}_0}} = 16}$
\cite[Theorem~6.6]{EJ14}. As~the geometric Picard rank does not increase under generisation, the claim is~established.\smallskip

\noindent
ii)
Spreading~out, one finds a morphism
$\smash{\underline{q}\colon \underline{X} \to \underline{B}}$
of schemes of finite type
over~$\bbZ$,
where
$\smash{\underline{B} \subset \Ab^1_{\bbZ[\frac1l]}}$
is an open subscheme. Restricting~to a further open subscheme, one may assume that every fibre is a
$K3$~surface.

Let~us put
$D := 8$
and
$a := 3$
or~$5$
and consider an arbitrary
\mbox{$\smash{\bbF_{\!p}}$-rational}
point
$\smash{\tau\in\underline{B}(\bbF_{\!p})}$,
for
$p \equiv a \pmod 8$
any prime~number. Then~one has
$\smash{\sqrt{2} \not\in \bbF_{\!p}}$.
From~this, one directly deduces that Frobenius operates on the six lines of the ramification locus by a permutation of
type~$(12)(34)(5)(6)$.
Therefore,~the induced operation on the fifteen pairs fixes only three of them,
$\{1,2\}$,
$\{3,4\}$,
and~$\{5,6\}$,
while the other twelve form six
$2$-cycles.

On~the other hand, write
$\underline{X}'_\tau$
for the double cover
of~$\Pb^2$
underlying~$\underline{X}_\tau$.
Then,~since
$p \equiv 3, 5 \pmod 8$,
it is known \cite[Theorem~6.3]{EJ14} that
$\#\underline{X}'_\tau(\bbF_{\!p}) = p^2 + p + 1$.
Furthermore, the lines
$E_{12}, \ldots, E_{56}$
are blown down, exactly three of which contain
\mbox{$\bbF_{\!p}$-rational}
points. Consequently,
$$\#\underline{X}_\tau(\bbF_{\!p}) = \#\underline{X}'_\tau(\bbF_{\!p}) + 3p = p^2 + 4p + 1 \equiv 1 \pmod p \,.$$
In~other words, Theorem~\ref{suff_crit}.a) applies and shows that the generic fibre
$X_\eta$
has indeed real or complex~multiplication.
Write~$E$
for the endomorphism field
of~$X_\eta$.

Moreover,~for
$t_0 \equiv 1 \pmod {17\!\cdot\!23}$,
the special fibre
$X_{t_0}$
has real multiplication by
$\smash{\bbQ(\sqrt{2})}$,
according to~\cite[Theorem~6.6]{EJ14}. Thus,~Corollary~\ref{endfield_spec} yields
$\smash{E \subseteq \bbQ(\sqrt{2})}$.
Together~with the fact that
$\smash{E \supsetneqq \bbQ}$,
this shows that
$X_\eta$
has real multiplication exactly
by~$\smash{\bbQ(\sqrt{2})}$.\smallskip

\noindent
iii)
This~is just an application of Theorem~\ref{suff_crit}.b).\smallskip

\noindent
iii)
Here,~the transcendental part
$\smash{T \subset H^2(X_\theta(\bbC), \bbQ)}$
is of
dimension~$6$.
By~Theorem \ref{suff_crit}.b), once again, about its endomorphism
field~$E_\theta$
we know that
$\smash{E_\theta \supseteq \bbQ(\sqrt{2})}$.

Assume~that equality does not~hold.
Then~$2 \mid [E_\theta\!:\!\bbQ] \mid 6$
and
$2 \neq [E_\theta\!:\!\bbQ]$,
which together leave
$[E_\theta\!:\!\bbQ] = 6$
as the only~option. In~this case, real multiplication is impossible, due to~\cite[Remark~1.5.3.c)]{Za}.
Thus,~$E_\theta$
must be a CM field of
degree~$6$.
Its~totally real subfield is hence cubic and contains
$\smash{\bbQ(\sqrt{2})}$,
a contradiction.%
\eop
\end{ex}

\section{A second explicit family}
\label{sec_exfam_Qw5}

In~order to specify the endomorphism field
$E \supsetneqq \bbQ$,
whose existence follows from Theorem~\ref{suff_crit}.a), a few particular assumptions are~necessary. We~will use the result below
for~$d=2$.

\begin{prop}
\label{estimates}
Let\/~$d$
and\/~$n$
be positive integers, 
$k$
a number field,
$X$
a\/
$K3$~surface
over\/~$k$,
and\/
$\frakp_1$
and\/~$\frakp_2$
be two prime ideals
of\/~$k$
at which\/
$X$
has good~reduction. Suppose~that

\begin{iii}
\item[$\bullet$ ]
$\rk\Pic X_{\overline{k}} \geq n$,
and
\item[$\bullet$ ]
$\smash{\rk\Pic X_{\overline\bbF_{\!\frakp_1}} = \rk\Pic X_{\overline\bbF_{\!\frakp_2}} = n+d}$,
whereas\/
$\smash{\disc \Pic X_{\overline\bbF_{\!\frakp_1}} / \disc \Pic X_{\overline\bbF_{\!\frakp_2}} \not\in (\bbQ^*)^2}$.
\end{iii}

\noindent
I.e., that the Picard lattices are incompatible in that sense that the quotient of their discriminants is a non-square
in\/~$\bbQ^*$.

\begin{abc}
\item
{\rm (Degree bound for the endomorphism field)}
Then, for the degree of the endomorphism
field\/~$E$,
the inequality\/
$[E\!:\!\bbQ] \leq d$
is~true.
\item
{\rm (Van Luijk's method)}
One has
$n \leq \rk\Pic X_{\overline{k}} \leq n+d-[E\!:\!\bbQ]$.
\end{abc}\smallskip

\noindent
{\bf Proof.}
{\em
a)
Applying van Luijk's method (cf.~\cite[Remark~3.2]{vL}) in the most naive way, one sees that
$\rk\Pic X_{\overline{k}} \leq n+d-1$.
Therefore,
\begin{equation}
\label{rk_ineq}
\smash{1 = (n+d)-(n+d-1) \leq \rk\Pic X_{\overline\bbF_{\!\frakp_1}}\! - \rk\Pic X_{\overline{k}} \leq (n+d)-n = d \,.}
\end{equation}
Lemma~\ref{rk_jump_mult} below shows that this yields
$[E\!:\!\bbQ] \leq d$.\smallskip

\noindent
b)
Again according to Lemma~\ref{rk_jump_mult}, one has that
$\smash{[E\!:\!\bbQ] \mid (\rk\Pic X_{\overline\bbF_{\!\frakp_1}}\! - \rk\Pic X_{\overline{k}})}$.
Thus, inequality~(\ref{rk_ineq}) shows that the difference between the two ranks is at least
$[E\!:\!\bbQ]$.
This~is exactly the asserted inequality to the~right. The inequality to the left is part of the~assumptions.
}
\eop
\end{prop}

\begin{lem}
\label{rk_jump_mult}
Let\/~$k$
be a number field,
$X$
a\/
$K3$~surface
over\/~$k$,
and\/
$\frakp$
be a prime ideal
of\/~$k$
at which\/
$X$
has good~reduction. Suppose~that\/
$X$
has real or complex multiplication by an endomorphism
field\/~$E$.
Then
$$\smash{[E\!:\!\bbQ] \mid (\rk\Pic X_{\overline\bbF_{\!\frakp}}\! - \rk\Pic X_{\overline{k}}) \,.}$$
{\bf Proof.}
{\em
Let~$\smash{\calP_\frakp \subseteq H^2_\et(X_{\overline{k}}, \bbQ_l(1))}$
be the vector space of those classes that are algebraic after specialisation
to~$\smash{H^2_\et(X_{\overline\bbF_{\!\frakp}}, \bbQ_l(1))}$,
and put
$\smash{\calT_\frakp := (\calP_\frakp)^\perp}$.
In~particular,
$\calT_\frakp \subseteq \calT$,
as algebraic classes remain algebraic under specialisation. To~specify this inclusion more precisely, let us write down an orthogonal decomposition
\begin{equation}
\label{splitting}
\calT = \calT_\frakp \oplus \frakS \,.
\end{equation}
I.e.,
$\frakS$
consists of the classes
in~$\calT$
that become algebraic after specialisation
to~$\smash{X_{\overline\bbF_{\!\frakp}}}$.
As~all these are defined over a finite extension
$\smash{\bbF_{\!\frakp^m}}$
of~$\smash{\bbF_{\!\frakp}}$,
for a certain
$m \in \bbN$,
the power
$\smash{(\Frob_\frakp)^m}$
operates
on~$\frakS$
with only
eigen\-value~$1$.

On~the other hand, there is
some~$n\in\bbN$
such that the operation of
$\smash{(\Frob_\frakp)^n}$
on
$\smash{\calT \subset H^2_\et(X_{\overline\bbF_{\!\frakp}}\!, \bbQ_l(1))}$
commutes with that
of~$E$,
cf.\ \cite[Corollary~4.2]{EJ14}. The~further power
$\smash{(\Frob_\frakp)^{nm}}$
operates
on~$\frakS$
again with only
eigen\-value~$1$.
However,~the Tate conjecture in the variant of Fact~\ref{Tate}.b) implies that
eigenvalue~$1$
does not occur
on~$\calT_\frakp$.
Since~$E$
and~$\smash{(\Frob_\frakp)^{nn'}}$
commute, this enforces that
$E$~maps
$\frakS$
to~itself.

Finally,~by Sublemma~\ref{descent_spec}, the splitting~(\ref{splitting}) descends
to~$T \subset H^2(X(\bbC), \bbQ)$,
$$T = T_\frakp \oplus S \,.$$
Moreover,~the operation
of~$E$
on~$T$
must map
$S$
to itself, as this is true after tensoring with
$\bbQ_l(1)$,
and
$\bbQ_l(1)$
is faithfully flat
over~$\bbQ$.
In~other words,
$S$~is
not just a
\mbox{$\bbQ$-vector}
space, but an
\mbox{$E$-vector}
space. Consequently,
$$[E\!:\!\bbQ] \mid \dim_\bbQ S = \rk\Pic X_{\overline\bbF_{\!\frakp}}\! - \rk\Pic X_{\overline{k}} \,,$$
as~claimed.
}
\eop
\end{lem}

\begin{subl}
\label{descent_spec}
Let\/~$k$
be a number field,
$X$
a\/
$K3$~surface
over\/~$k$,
and\/
$\frakp$
be a prime ideal
of\/~$k$
at which\/
$X$
has good~reduction. Denote
by\/~$\smash{\calP_\frakp \subseteq H^2_\et(X_{\overline{k}}, \bbQ_l(1))}$
the vector space of the classes being algebraic after specialisation
to\/~$\smash{H^2_\et(X_{\overline\bbF_{\!\frakp}}, \bbQ_l(1))}$.\smallskip

\noindent
Then~there is some subvector space\/
$\smash{P_\frakp \subset H^2(X(\bbC), \bbQ)}$
such that, under the comparison isomorphism\/ \cite[Expos\'e~11, Th\'eor\`eme~4.4.iii)]{SGA4},
$\calP_\frakp = P_\frakp \!\otimes_\bbQ\! \bbQ_l(1)$.\medskip

\noindent
{\bf Proof.}
{\em
For
$A$
any
\mbox{$\bbQ_l$-vector}
space, the Chern class homomorphism, combined with the specialisation and comparison isomorphisms, yields an~injection
$$c_1^{(A)}\colon \Pic X_{\overline\bbF_{\!\frakp}} \!\otimes_\bbZ\! A \longrightarrow H^2_\et(X_{\overline\bbF_{\!\frakp}}, A) \cong H^2_\et(X_{\overline{k}}, A) \cong H^2(X(\bbC), \bbQ) \!\otimes_\bbQ\! A \,,$$
which is natural
in~$A$.
In~particular, the difference kernels of the vertical arrows of the commutative~diagram
$$
\xymatrixrowsep{3mm}
\xymatrix{
\Pic X_{\overline\bbF_{\!\frakp}} \!\otimes_\bbZ\! \bbQ_l(1) \!\otimes_\bbQ\! \bbQ_l(1) \ar[rr]^{\!\!\!\!\!\!\!\!\!\!c_1^{(\bbQ_l(1) \otimes_\bbQ\! \bbQ_l(1))}} && H^2(X(\bbC), \bbQ) \!\otimes_\bbQ\! \bbQ_l(1) \!\otimes_\bbQ\! \bbQ_l(1) \\
\Pic X_{\overline\bbF_{\!\frakp}} \!\otimes_\bbZ\! \bbQ_l(1) \ar@<-.5ex>[u]^{\id\!\otimes\!1\!\otimes\id\;\;} \ar@<.5ex>[u]_{\;\;\id\!\otimes\!\id\!\otimes1} \ar[rr]^{c_1^{(\bbQ_l(1))}} && H^2(X(\bbC), \bbQ) \!\otimes_\bbQ\! \bbQ_l(1) \ar@<-.5ex>[u]^{\id\!\otimes\!1\!\otimes\id\;\;} \ar@<.5ex>[u]_{\;\;\id\!\otimes\!\id\!\otimes1}
}
$$
are connected by a descent homomorphism
$$c_\frakp\colon \Pic X_{\overline\bbF_{\!\frakp}} \!\otimes_\bbZ\! \bbQ \longrightarrow H^2(X(\bbC), \bbQ) \,.$$
Its~image
$\im c =: P_\frakp$
is the desired subvector~space.
}
\eop
\end{subl}

\begin{nota}
Let\/~$X$
be a\/
$K3$~surface
over a number
field\/~$k$
and\/
$\frakp$
a prime ideal
of\/~$k$,
at which\/
$X$
has good reduction.

\begin{iii}
\item
We~denote by
$\smash{\chi^\calT_{\frakp^n}}$
the characteristic polynomial of
$(\Frob_\frakp)^n \in \Gal(\overline{k}/k)$
on the transcendental
part~$\calT \subset H^2_\et(X_{\overline{k}}, \bbQ_l(1))$.
\item
We~factorise
$\smash{\chi^\calT_{\frakp^n} \in \bbQ[Z]}$
completely in the form
$\chi^\calT_{\frakp^n}(Z) = \chi^\tr_{\frakp^n}(Z) \cdot \prod_{i=1}^d (Z - \zeta_{k_i}^{e_i})$,
for
$k_1,\ldots,k_d \in \bbN$.
I.e.\ in such a way that
$\smash{\chi^\tr_{\frakp^n} \in \bbQ[Z]}$
does not have any further zeroes being roots of~unity.
\end{iii}
\end{nota}

\begin{rems}
\begin{iii}
\item
The power
$(\Frob_\frakp)^n$
of the Frobenius element is uniquely determined up to conjugation, so
$\smash{\chi^\calT_{\frakp^n}}$
is well-defined. By~\cite[Expos\'e~XVI, Corollaire~2.2]{SGA4}, it coincides with the characteristic polynomial of
$\Frob$
on the corresponding part
of~$\smash{H^2_\et(X_{\overline\bbF_{\!\frakp}}, \bbQ_l(1))}$.
In~particular,
$\smash{\chi^\calT_{\frakp^n}}$
is, in fact, a polynomial with coefficients
in~$\bbQ$
\cite[Th\'eor\`eme~1.6]{De74}.

By~definition, one has
$\smash{\deg\chi^\calT_{\frakp^n} = 22 - \rk\Pic X_{\overline{k}}}$.
\item
According~to the Tate conjecture (Fact~\ref{Tate}.b)),
$\smash{\chi^\tr_{\frakp^n}}$
is the characteristic polynomial of
$\Frob^n$
on the transcendental part of
$\smash{H^2_\et(X_{\overline\bbF_{\!\frakp}}\!, \bbQ_l(1))}$.
Let~us note, in~particular, that
$\smash{\deg\chi^\tr_{\frakp^n} = 22 - \rk\Pic X_{\overline\bbF_{\!\frakp}}}$.
\item
Consequently,~one
has that
$d = \rk\Pic X_{\overline\bbF_{\!\frakp}}\! - \rk\Pic X_{\overline{k}}$.
\end{iii}
\end{rems}

In~order to decide which quadratic number field exactly is the endomorphism field, the following result is~useful.

\begin{prop}
\label{factorisation_E}
Let\/~$\frakp$
be a prime of good reduction of the\/
$K3$~surface\/~$X$
over a number field\/~$k$,
having real or complex multiplication by a
field\/~$E$
containing the quadratic number
field\/~$\smash{\bbQ(\sqrt{\delta})}$.
Then~at least one of the following two statements is~true.

\begin{iii}
\item
The polynomial\/
$\chi^\tr_\frakp \in \bbQ[Z]$
is the norm of a polynomial from\/
$\smash{\bbQ(\sqrt{\delta})[Z]}$.
\item
For~some\/
$n\in\bbN$,
the polynomial\/
$\smash{\chi^\tr_{\frakp^n}}$
is a square
in\/~$\bbQ[Z]$.
\end{iii}\smallskip

\noindent
{\bf Proof.}
{\em
This is essentially \cite[Theorem~4.9]{EJ14}. Note that the proof given in \cite{EJ14} works over an arbitrary number~field.
}
\eop
\end{prop}

\begin{nota}
Let
$\smash{q\colon X \to B}$
be the family from Example~\ref{Qw5_family} and
$\smash{q'\colon X' \to B}$
the underlying family of double covers
of~$\Pb^2$.

\begin{iii}
\item
Spread~out in the obvious way, i.e.\ put\/
$\smash{\underline{B} := \Spec\bbZ[T\!,\!\frac1{(T-1)(T^4-T^3+T^2-T+1)}] \!\subset\! \Ab^1_\bbZ}$
and let\/
$\smash{\underline{X}'}$
be the double cover of\/
$\smash{\Pb^2_{\underline{B}}}$
given by the same equation as~(\ref{Qw5}). We~denote the family of schemes thus obtained
by~$\smash{\underline{q}'\colon \underline{X}' \to \underline{B}}$.
\item
Let~$\smash{\underline{X}}$
be the blow-up of
$\smash{\underline{X}'}$
in the 15
\mbox{codimension-$2$}
subschemes being the Zariski closures of the 15 singular points on the generic fibre~$\smash{\underline{X}'_\eta}$.
For~the resulting family, let us write
$\smash{\underline{q}\colon \underline{X} \to \underline{B}}$.
\end{iii}
\end{nota}

\begin{proo3}
For~$t_0 \in B$,
the ramification locus of the double cover
underlying~$X_{t_0}$
is a union of six lines, all of which are defined over
$\smash{\bbQ(\zeta_5,t_0)}$.
Indeed,~the quartic occurring in (\ref{Qw5}) is the norm form of the linear
form~$x - \zeta_5y + \zeta_5^2z$
defined
over~$\smash{\bbQ(\zeta_5)}$.

Moreover,~no three of the five lines in the ramification locus that do not depend
of~$t$
have a point in common, the ten points of intersection being
$\smash{(1\!:\! \frac{\sqrt{5}-1}2 \!:\! 1)}$,
$(1 \!:\! \zeta_5 \!+\! \zeta_5^3 \!:\! \zeta_5^4)$,
$(1\!:\!0\!:\!-\zeta_5)$,
and their~conjugates. A~direct calculation shows that the sixth line passes through one of these points if and only of
$t_0=1$
or
$t_0$~is
a proper tenth root of unity. As~these are exactly the points excluded
from~$B$,
$q$~is
indeed a family of
$K3$~surfaces.\smallskip

\noindent
i)
Having~spread out as described, one has the morphism
$\smash{\underline{q}\colon \underline{X} \to \underline{B}}$
of schemes of finite type
over~$\bbZ$.
Restricting~to a suitable open subscheme, one may assume that
$\smash{\underline{B} \subseteq \Ab^1_{\bbZ[\frac1l]}}$
and that every fibre is a
$K3$~surface.

Furthermore, since
$\Pic X_{\overline\eta}$,
contains the span of the
$\smash{(\atop{6}2) = 15}$
classes of the exceptional curves, obtained by blowing up the intersection points of the ramification locus, and the pull-back of the class of the general line
on~$\Pb^2$,
one certainly has
$\rk\Pic X_{\overline\eta} \geq 16$.

In~order to apply Theorem~\ref{suff_crit}, let us put
$D := 5$
and
$a := 2$
or~$3$
and consider an arbitrary
\mbox{$\bbF_{\!p}$-rational}
point
$\tau\in\underline{B}(\bbF_{\!p})$,
for
$p \equiv a \pmod 5$
any prime number.
Then~$\smash{\sqrt{5} \not\in \bbF_{\!p}}$.
Therefore,~Frobenius operates on the six lines of the ramification locus by a permutation of
type~$(1)(2)(3456)$.
Hence,~the induced operation on the fifteen pairs fixes only one of them,
$\{1,2\}$,
while the other twelve form two
$4$-cycles
and two
$2$-cycles.

On~the other hand, Lemma~\ref{Qw5_count} below shows that
$\#\underline{X}'_\tau(\bbF_{\!p}) = p^2 + p + 1$.
As~the lines
$E_{12}, \ldots, E_{56}$
are blown down, exactly one of which contains
\mbox{$\bbF_{\!p}$-rational}
points, this~yields
$$\#\underline{X}_\tau(\bbF_{\!p}) = \#\underline{X}'_\tau(\bbF_{\!p}) + p = p^2 + 2p + 1 \,.$$
In~other words, Theorem~\ref{suff_crit}.a) applies and shows that
$X_\eta$
has indeed real or complex~multiplication.
Write~$E$
for the endomorphism field
of~$X_\eta$.
Clearly,~$[E\!:\!\bbQ] \geq 2$.

\looseness-1
Next,~consider the closed point
$t_0 := 15 \in B$.
Then~$\smash{\rk\Pic X_{\overline{t}_0} \geq \rk\Pic X_{\overline\eta} \geq 16}$,
as the geometric Picard rank cannot drop under specialisation. Moreover,~by Corollary~\ref{endfield_spec}, for the endomorphism field, one has
$E_{t_0} \supseteq E$
and, in particular,
$[E_{t_0}\!:\!\bbQ] \geq 2$.
On~the other hand,
$X_{t_0}$
has two reductions of geometric Picard
rank~$18$
with incompatible discriminants, cf.\ Lemma \ref{Qw5_reds}. Thus,~Proposition~\ref{estimates} applies
to~$X_{t_0}$
with
$n=16$
and~$d=2$.
It~yields that
$\smash{\rk\Pic X_{\overline{t}_0} = 16}$
and that
$[E_{t_0}\!:\!\bbQ] = 2$.

As~$\eta$
is a generisation
of~$t_0$,
one finds that
$\smash{\rk\Pic X_{\overline\eta} \leq 16}$.
As~we saw the other inequality above, the proof of~i) is~complete.\smallskip

\noindent
ii)
Recall~the facts that
$[E\!:\!\bbQ] \geq 2$,
$[E_{t_0}\!:\!\bbQ] = 2$,
and
$E \subseteq E_{t_0}$,
which were all found during the proof of~i). Together,~they immediately show that
$[E\!:\!\bbQ] = 2$.
Consequently,
$E_{t_0} = E$,
so it is enough to verify
that~$\smash{E_{t_0} = \bbQ(\sqrt{5})}$.

As~$E_{t_0}$
is known to be a quadratic number field, this follows from Proposition~\ref{factorisation_E}. Indeed,
$$\textstyle \chi^\tr_{19} = Z^4 - \frac{14}{19}Z^3 + \frac{34}{19}Z^2 - \frac{14}{19}Z + 1 = (Z^2 - \frac{7+5\sqrt{5}}{19}Z + 1)(Z^2 - \frac{7-5\sqrt{5}}{19}Z + 1)$$
splits
over~$\smash{\bbQ(\sqrt{5})}$
and over no other quadratic~field. In~fact,
$\smash{\Gal(\chi^\tr_{19}) = D_4}$,
which has only one intransitive subgroup of index~two. Observe,~moreover, that the splitting field of
$\smash{\chi^\tr_{19}}$
does not contain any roots of unity, except
for~$(-1)$.
Thus,~if
$\smash{\chi^\tr_{19^n}}$
were a perfect square for some
$n\in\bbN$
then this would happen for
$n=2$,
already, which is not the~case.\smallskip

\noindent
iii)
According~to Corollary~\ref{complex}, this is a direct consequence of~ii).\smallskip

\noindent
iv)
This is exactly the same argument as in Example~\ref{Qw2_family}.iv).
\eop
\end{proo3}

\begin{lem}
\label{Qw5_reds}
\begin{abc}
\item
Put\/
$\tau_1 := (0 \bmod 3) \in \underline{B} \subset \Ab^1_\bbZ$.
Then~the special fibre\/
$\underline{X}_{\tau_1}$
of\/~$\smash{\underline{q}}$
is a\/
$K3$~surface
over\/~$\bbF_{\!3}$
of geometric Picard
rank\/~$18$.
The~discriminant of the Picard lattice
is\/~$\smash{\overline{(-1)} \in \bbQ^*/(\bbQ^*)^2}$
and one has\/
$\chi^\tr_3(Z) = Z^4 - \frac43Z^2 + 1$.
\item
Put\/
$\tau_2 := (15 \bmod 19) \in \underline{B} \subset \Ab^1_\bbZ$.
Then~the special fibre\/
$\underline{X}_{\tau_2}$
of\/~$\smash{\underline{q}}$
is a\/
$K3$~surface
over\/~$\bbF_{\!19}$
of geometric Picard
rank\/~$18$.
The discriminant of the Picard lattice
is\/~$\smash{\overline{(-11)} \in \bbQ^*/(\bbQ^*)^2}$
and one has\/
$\chi^\tr_{19}(Z) = Z^4 - \frac{14}{19}Z^3 + \frac{34}{19}Z^2 - \frac{14}{19}Z + 1$.
\end{abc}\smallskip

\noindent
{\bf Proof.}
{\em
One~first has to check that the
family~$\smash{\underline{q}}$
has good reduction at these two closed~points. For~this, let us note the following.

The five lines in the ramification locus that are independent
of~$t$
are distinct in every
characteristic~$\neq \!5$.
Their~ten points of intersection are distinct in every
characteristic~$\neq \!2, 5$.
Moreover,~the sixth ramification line is distinct from the others, for every value
of~$t$,
in every
characteristic~$\neq \!5$.
Thus,~if
$p \neq 2,5$,\,
$a \not\equiv 1 \pmod p$,
and
$a^5 \not\equiv -1 \pmod p$
then
$q$
has good reduction at
$(a \bmod p) \in \underline{B} \subset \Ab^1_\bbZ$.
This~criterion applies to
$\tau_1 = (0 \bmod 3)$,
as well as to
$\tau_2 = (15 \bmod 19)$.

In~order to determine the characteristic polynomials and geometric Picard ranks, one has to count the points on
$\underline{X}_{\tau_1}$
and~$\underline{X}_{\tau_2}$
that are defined over the prime field and some of its extensions.
We~applied routine methods for this, using {\tt magma}. For~some background, the reader might consult~\cite{EJ16}. Finally,~the discriminants are easily calculated, using the Artin-Tate formula~\cite[Theorem~6.1]{Mi}.
}
\eop
\end{lem}

\begin{lem}
\label{Qw5_count}
For~every prime number\/
$p \equiv 2,3 \pmod 5$
and every\/
\mbox{$\bbF_{\!p}$-rational}
point\/
$\smash{\tau\in \underline{B}(\bbF_{\!p})}$,
the special fibre\/
$\smash{\underline{X}'_\tau}$
of\/~$\smash{\underline{q}'}$
has point~count
$$\#\underline{X}'_\tau(\bbF_{\!p}) = p^2 + p + 1 \,.$$
{\bf Proof.}
{\em
Write
$\smash{\widetilde\Pb{}^2_{\bbF_{\!p}}}$
for the blow-up
of~$\smash{\Pb^2_{\bbF_{\!p}}}$
in~$(2:(-1):2)$
and put
$\smash{\widetilde{X} := \underline{X}' \times_{\Pb^2_{\underline{B}}} \widetilde\Pb{}^2_{\bbF_{\!p}}}$.
Since~$(2:(-1):2)$
is, independently of the value
of~$\tau$,
a point on the ramification locus, the assertion is equivalent to
$\smash{\#\widetilde{X}_\tau(\bbF_{\!p}) = p^2 + 2p + 1}$.

Moreover,
$\smash{\widetilde\Pb{}^2_{\bbF_{\!p}}}$
is fibred into the lines through
$(2\!:\!(-1)\!:\!2)$.
We~parametrise the fibration
by~$(v_0\!:\!v_1)\in\Pb^1$
and let
$l_{(v_0:v_1)}$
be the line parametrised~by
$$u \mapsto ((-2u\!+\!v_0) \!:\! (u\!+\!v_1) \!:\! (-2u)) \,.$$
Note~here that this line does not depend on the choice of representatives
for~$(v_0\!:\!v_1)$
and that all these lines have
$(2\!:\!(-1)\!:\!2)$
as their point at~infinity.

Correspondingly,
$\smash{\widetilde{X}_\tau}$
is fibred into
\mbox{genus-$2$-curves}
$\smash{C_{\tau,(v_0:v_1)}}$.
It~is clearly sufficient to show that
$\#C_{\tau,(v_0:v_1)}(\bbF_{\!p}) = p+1$,
for every
$\tau \in \underline{B}(\bbF_{\!p})$
and every
$(v_0\!:\!v_1) \in \Pb^1(\bbF_{\!p})$.

A~direct calculation shows that these curves are given~by
$$C_{\tau,(v_0:v_1)}\colon~ w^2 = 25 (v_0\!+\!(-2t\!+\!2)v_1) \!\cdot\! P_{(v_0,v_1)}(u)$$
for
\begin{align}
\label{ram_normed}
P_{(v_0,v_1)}(u) := u^5 \!+\! (-v_0 &\!+\! v_1)u^4 \\[-1mm]
 &\!+\! \textstyle (\frac35v_0^2 \!-\! \frac65v_0v_1 \!-\! \frac25v_1^2)u^3 \!+\! (-\frac15v_0^3 \!+\! \frac35v_0^2v_1 \!-\! \frac25v_1^3)u^2 \nonumber \\[-1mm]
 &\textstyle\hspace{1.5cm} + (\frac1{25}v_0^4 \!-\! \frac4{25}v_0^3v_1 \!+\! \frac1{25}v_0^2v_1^2 \!+\! \frac6{25}v_0v_1^3 \!+\! \frac1{25}v_1^4)u \nonumber \\[-1mm]
 &\textstyle\hspace{2.0cm} + \frac1{25}v_0^4v_1 \!+\! \frac1{25}v_0^3v_1^2 \!+\! \frac1{25}v_0^2v_1^3 \!+\! \frac1{25}v_0v_1^4 \!+\! \frac1{25}v_1^5 \,. \nonumber
\end{align}
Here, the coefficient
$v_0\!+\!(-2t\!+\!2)v_1$
might be~zero. Then~the curve is a double line and
$\#C_{\tau,(v_0:v_1)}(\bbF_{\!p}) = p+1$
is obviously~true. Otherwise,
$C_{\tau,(v_0:v_1)}$
is a quadratic twist of the
curve~$C_{(v_0:v_1)}$,
given~by
$$C_{(v_0:v_1)}\colon~ w^2 = P_{(v_0,v_1)}(u) \,,$$
so~it suffices to show that
$\#C_{(v_0:v_1)}(\bbF_{\!p}) = p+1$.

For~this, we note that
$P_{(v_0,v_1)}$
is a permutation polynomial, according to Sublemma~\ref{perm_pol}. Thus~the number of points
on~$C_{(v_0:v_1)}$
is the same as that on the curve, given by
$w^2 = u$,
which is isomorphic to the projective~line.
}
\eop
\end{lem}

\begin{subl}
\label{perm_pol}
Let\/~$p \equiv 2,3 \pmod 5$
be a prime~number. Then,~for arbitrary\/
$v_0, v_1 \in \bbF_{\!p}$,
the quintic polynomial\/
$P_{(v_0,v_1)} \in \bbF_{\!p}[u]$,
defined in~(\ref{ram_normed}), is a permutation polynomial. I.e.,~it
induces a bijection of\/
$\bbF_{\!p}$
onto~itself.\medskip

\noindent
{\bf Proof.}
{\em
Define
$\smash{\widetilde{P} \in \bbF_{\!p}[u]}$
by
$\smash{\widetilde{P}(u) := P_{(v_0,v_1)}(u+\frac{v_0-v_1}5) - C}$,
for
$\smash{C := P_{(v_0,v_1)}(\frac{v_0-v_1}5)}$.
This~is a normalised form
of~$P$,
cf.~\cite[the remarks before Theorem~7.11]{LN}. It~is clearly sufficient to show that
$\smash{\widetilde{P}}$
is a permutation poly\-no\-mial.

For~this, a direct calculation shows that
\begin{align*}
\widetilde{P}(u) &\textstyle= u^5 + (\frac15v_0^2 \!-\! \frac25v_0v_1 \!-\! \frac45v_1^2)u^3 + (\frac1{125}v_0^4 \!-\! \frac4{125}v_0^3v_1 \!-\! \frac4{125}v_0^2v_1^2 \!+\! \frac{16}{125}v_0v_1^3 \!+\! \frac{16}{125}v_1^4)u \\
&= u^5 - 5\alpha u^3 + 5\alpha^2 u \,,
\end{align*}
when putting
$\alpha := -\frac1{25}v_0^2 + \frac2{25}v_0v_1 + \frac4{25}v_1^2$.
This~means that
$\smash{\widetilde{P}}$
coincides with the Dickson polynomial
$g_5(u, \alpha)$,
which is known \cite[Theorem~7.16]{LN} to be a permutation polynomial, for every prime
number~$p \equiv 2,3 \pmod 5$.
The proof is therefore~complete.%
}%
\eop
\end{subl}

\setlength\parindent{0mm}
\end{document}